\documentclass{amsart}
\usepackage{amssymb,amsmath,comment,mathrsfs,mathtools,tikz-cd,todonotes}
\usepackage[colorlinks=true,linkcolor=blue]{hyperref}
\usepackage{cleveref}

\theoremstyle{remark}

\newtheorem{para}{\bf}[subsection]

\newtheorem{example}[para]{\bf Example}

\newtheorem{rem}[para]{\bf Remark}

\theoremstyle{definition}

\newtheorem{dfn}[para]{Definition}

\theoremstyle{plain}

\newtheorem{theorem}[para]{Theorem}

\newtheorem{lemma}[para]{Lemma}
\newtheorem{cor}[para]{Corollary}

\newtheorem{prop}[para]{Proposition}

\newenvironment{numequation}{\addtocounter{para}{1}
\begin{equation}}{\end{equation}}

\newcommand{\bbP}{{\mathbb P}}
\newcommand{\bbQ}{{\mathbb Q}}
\newcommand{\bbR}{{\mathbb R}}

\newcommand{\bbZ}{{\mathbb Z}}

\newcommand{\bB}{{\bf B}}

\newcommand{\bG}{{\bf G}}

\newcommand{\bP}{{\bf P}}

\newcommand{\bT}{{\bf T}}

\newcommand{\bZ}{{\bf Z}}

\newcommand{\frb}{{\mathfrak b}}

\newcommand{\frg}{{\mathfrak g}}

\newcommand{\frl}{{\mathfrak l}}
\newcommand{\frm}{{\mathfrak m}}
\newcommand{\frn}{{\mathfrak n}}

\newcommand{\frp}{{\mathfrak p}}

\newcommand{\frs}{{\mathfrak s}}
\newcommand{\frt}{{\mathfrak t}}

\newcommand{\frx}{{\mathfrak x}}

\newcommand{\frz}{{\mathfrak z}}

\newcommand{\cA}{{\mathcal A}}
\newcommand{\cB}{{\mathcal B}}

\newcommand{\cF}{{\mathcal F}}

\newcommand{\cO}{{\mathcal O}}

\newcommand{\cT}{{\mathcal T}}

\newcommand{\sE}{{\mathscr E}}
\newcommand{\sF}{{\mathscr F}}
\newcommand{\sG}{{\mathscr G}}

\newcommand{\sL}{{\mathscr L}}

\newcommand{\sS}{{\mathscr S}}
\newcommand{\sT}{{\mathscr T}}

\newcommand{\sZ}{{\mathscr Z}}

\newcommand{\pair}[2]{\langle #1 , #2 \rangle}

\newcommand{\Ad}{{\rm Ad}}

\newcommand{\alg}{{\rm alg}}

\newcommand{\bksl}{\backslash}

\newcommand{\diag}{{\rm diag}}

\newcommand{\End}{{\rm End}}

\newcommand{\Ex}{{E^\times}}
\newcommand{\Ext}{{\rm Ext}}

\newcommand{\Fx}{{F^\times}}

\newcommand{\GL}{{\rm GL}}

\newcommand{\Hom}{{\rm Hom}}

\newcommand{\hra}{\hookrightarrow}
\newcommand{\id}{{\rm id}}

\newcommand{\Ind}{{\rm Ind}}

\newcommand{\irr}{\rm irr}
\newcommand{\la}{{\rm la}}

\newcommand{\Lie}{{\rm Lie}}

\newcommand{\lra}{\longrightarrow}
\newcommand{\Max}{{\rm Max}}

\newcommand{\midc}{{\; | \;}}

\newcommand{\ot}{\otimes}

\newcommand{\pr}{{\rm pr}}

\newcommand{\Qp}{{\bbQ_p}}

\newcommand{\Qpx}{{\bbQ^\times_p}}
\newcommand{\ra}{\rightarrow}
\newcommand{\Rep}{{\rm Rep}}

\newcommand{\rmd}{{\rm d}}

\newcommand{\sm}{{\rm sm}}

\renewcommand{\sp}{{\rm sp}}

\newcommand{\sub}{\subset}

\newcommand{\triv}{{\bf 1}}

\newcommand{\vep}{{\varepsilon}}

\newcommand{\vpi}{{\varpi}}
\newcommand{\vtheta}{{\vartheta}}
\newcommand{\wdelta}{\widehat{\delta}}

\newcommand{\wlambda}{\widetilde{\lambda}}
\newcommand{\wmu}{\widetilde{\mu}}
\newcommand{\wnu}{\widetilde{\nu}}

\newcommand{\wot}{\widehat{\otimes}}

\newcommand{\x}{{\times}}

\newcommand{\Z}{{\mathbb Z}}

\newcommand{\st}{\hspace{0.1cm}|\hspace{0.1cm}}

\newcommand{\chcFGP}{{\check{\cF}^G_P}}

\newcommand{\cOp}{{\cO^\frp}}
\newcommand{\cOpinfty}{{\cO^{\frp,\infty}}}
\newcommand{\cOpinftyalg}{{\cO^{\frp,\infty}_\alg}}
\newcommand{\cOPinfty}{{\cO^{P,\infty}}}

\newcommand{\cTlamu}{\cT_\lambda^\mu}
\newcommand{\cTmula}{\cT_\mu^\lambda}
\newcommand{\cTclamu}{\cT^\mu_{c,\lambda}}
\newcommand{\cTcmula}{\cT^\lambda_{c,\mu}}
\newcommand{\cTrlamu}{\cT^\mu_{r,\lambda}}
\newcommand{\cTrmula}{\cT^\lambda_{r,\mu}}

\newcommand{\DB}{{D(B)}}

\newcommand{\DG}{{D(G)}}
\newcommand{\DGcoad}{{D(G)}\mbox{-}{\rm mod}^{\rm coad}}
\newcommand{\DGcoadzfin}{{D(G)}\mbox{-}{\rm mod}^{\rm coad}_{\frz{\mbox{-}{\rm fin}}}}
\newcommand{\DGcoadlambda}{{D(G)}\mbox{-}{\rm mod}^{\rm coad}_{|\lambda|}}
\newcommand{\DGcoadmu}{{D(G)}\mbox{-}{\rm mod}^{\rm coad}_{|\mu|}}
\newcommand{\DGmod}{{D(G)}\mbox{-}{\rm mod}}
\newcommand{\DgB}{D(\mathfrak{g},B)}
\newcommand{\DgP}{D(\mathfrak{g},P)}
\newcommand{\Dmu}{D(G)\mbox{-}{\rm mod}_{|\mu|}}

\renewcommand{\DH}{{D(H)}}

\newcommand{\Dlambda}{D(G)\mbox{-}{\rm mod}_{|\lambda|}}

\newcommand{\Dzfin}{D(G)\mbox{-}{\rm mod}_{\frz{\mbox{-}{\rm fin}}}}

\newcommand{\DP}{{D(P)}}
\newcommand{\DrH}{{D_r(H)}}
\newcommand{\DrHfglambda}{D_r(H)\mbox{-}{\rm mod}^{\rm fg}_{|\lambda|}}
\newcommand{\DrHfgmu}{D_r(H)\mbox{-{\rm mod}}^{\rm fg}_{|\mu|}}
\newcommand{\DrHzfin}{D_r(H)\mbox{-}{\rm mod}_{\frz{\mbox{-}{\rm fin}}}}
\newcommand{\DrHfgzfin}{D_r(H)\mbox{-}{\rm mod}^{\rm fg}_{\frz{\mbox{-}{\rm fin}}}}
\newcommand{\DrHmod}{D_r(H)\mbox{-}{\rm mod}}

\newcommand{\DrHlambda}{D_r(H)\mbox{-}{\rm mod}_{|\lambda|}}
\newcommand{\DrHmu}{D_r(H)\mbox{-}{\rm mod}_{|\mu|}}

\newcommand{\hsT}{\widehat{\sT}}
\newcommand{\Lift}{{\rm Lift}}

\newcommand{\onu}{\overline{\nu}}
\newcommand{\Phibrla}{\Phi_{[\lambda]}}
\newcommand{\Rmod}{R\mbox{-}{\rm mod}}
\newcommand{\Rzfin}{R\mbox{-}{\rm mod}_{\frz{\mbox{-}{\rm fin}}}}
\newcommand{\Rfrm}{R\mbox{-}{\rm mod}_\frm}
\newcommand{\Rlambda}{R\mbox{-}{\rm mod}_{|\lambda|}}
\newcommand{\sadm}{{{\rm s}\mbox{-}{\rm adm}}}

\newcommand{\Tlamu}{T_\lambda^\mu}
\newcommand{\Tmula}{T_\mu^\lambda}

\newcommand{\UbE}{{U(\frb_E)}}
\newcommand{\Ug}{{U(\frg)}}
\newcommand{\UgE}{{U(\frg_E)}}
\newcommand{\Ugzfin}{U(\frg_E)\mbox{-}{\rm mod}_{\frz{\mbox{-}{\rm fin}}}}
\newcommand{\Uglambda}{U(\frg_E)\mbox{-}{\rm mod}_{|\lambda|}}
\newcommand{\Ugmu}{U(\frg_E)\mbox{-}{\rm mod}_{|\mu|}}
\newcommand{\UgEmod}{U(\frg_E)\mbox{-{\rm mod}}}
\newcommand{\Ugfglambda}{U(\frg_E)\mbox{-}{\rm mod}^{{\rm fg}}_{|\lambda|}}
\newcommand{\Ugfgmu}{U(\frg_E)\mbox{-}{\rm mod}^{{\rm fg}}_{|\mu|}}
\newcommand{\Ugfrm}{U(\frg_E)\mbox{-}{\rm mod}_\frm}
\newcommand{\Ugmod}{U(\frg)\mbox{-{\rm mod}}}

\newcommand{\UpE}{{U(\frp_E)}}
\newcommand{\Wbrla}{{W_{[\lambda]}}}
\newcommand{\whs}{\widehat{s}}
\begin{document}

\title{Translation functors for locally analytic representations}

\author{Akash Jena}
\address{Indiana University, Department of Mathematics, Rawles Hall, Bloomington, IN 47405, U.S.A.}
\email{akjena@iu.edu}

\author{Aranya Lahiri}
\address{University of California at San Diego, Department of Mathematics, 9500 Gilman Dr, La Jolla, CA 92093, U.S.A.}
\email{arlahiri@ucsd.edu}

\author{Matthias Strauch}
\address{Indiana University, Department of Mathematics, Rawles Hall, Bloomington, IN 47405, U.S.A.}
\email{mstrauch@indiana.edu}

\begin{abstract}
Let $G$ be a $p$-adic Lie group with reductive Lie algebra $\frg$.  In analogy to the translation functors introduced by Bernstein and Gelfand on categories of $U(\frg)$-modules we consider similarly defined functors on the category of modules over the locally analytic distribution algebra $\DG$ on which the center of $U(\frg)$ acts locally finite. These functors induce equivalences between certain subcategories of the latter category. Furthermore, these translation functors are naturally related to those on category $\cO$ via the functors from category $\cO$ to the category of coadmissible modules. We also investigate the effect of the translation functors on locally analytic representations $\Pi(V)^\la$ associated by the $p$-adic Langlands correspondence for $\GL_2(\Qp)$ to 2-dimensional Galois representations $V$. 
\end{abstract}
\maketitle

\tableofcontents

\section{Introduction}
\setcounter{subsection}{1}
Translation functors in the sense of Jantzen \cite{JantzenModuln} or Bernstein-Gelfand \cite{BeGe} are an important tool for understanding the category of modules over a semisimple Lie algebra $\frg$. In \cite{JantzenModuln} the Lie algebra is assumed to be split over a field of characteristic zero, which we here denote by $F$ and which, for our purposes, will be a finite extension of $\Qp$. We also fix a coefficient field $E/F$ and consider representations of $\frg$ on $E$-vector spaces, or, equivalently, $\UgE$-modules.\footnote{The introduction of $E$ at this point seems superfluous, as we could take $E=F$ here, but when we consider $p$-adic groups later on, this setting provides more flexibility.} 

\vskip8pt

{\it Translation functors for $\frg$-modules.} We start by recalling the translation functors for categories of $\UgE$-modules. Let $\frz_E$ be the center of $\UgE$ and consider the category $\Ugzfin$ of modules $M$ such that every $m \in M$ is annihilated by an ideal of finite codimension in $\frz_E$. Modules $M$ in this category split into their $\frm$-primary components: $M = \bigoplus_{\frm \in \Max(\frz_E)} M_\frm$, and we denote by $\Ugfrm$ the subcategory of modules $M$ with the property that $M = M_\frm$. The projection $M \ra M_\frm$ is denoted by $\pr_\frm$. Translation functors are endo-functors of $\Ugzfin$ of the form

\[M \rightsquigarrow \pr_\frn (L \ot_E \pr_\frm(M)) \;,\]

\vskip8pt

where $L$ is a finite-dimensional irreducible representation of $\frg_E$, and $\frm,\frn \in \Max(\frz_E)$. 

\vskip8pt

Following the set-up in \cite{JantzenModuln} and \cite{Hu}, we usually consider only those maximal ideals $\frm \sub \frz_E$ with $\frz_E/\frm = E$. Fix a Borel subalgebra $\frb \sub \frt$ and an $F$-split Cartan sualgebra $\frt \sub \frb$. By the Harish-Chandra homomorphism (extended to the case of reductive $\frg$, cf. \cite[4.115]{KnappVogan_Cohomological}), any such $\frm$ is the kernel of the character $\chi_\lambda$ by which $\frz_E$ acts on the Verma module $M(\lambda) = \UgE \ot_{\UbE} E_\lambda$. As in the semisimple case one has $\chi_\lambda = \chi_\mu$ if and only if $\mu$ belongs the the orbit $|\lambda|$ of the Weyl group $W$ under the `dot-action'. We write $\pr_{|\lambda|}$ and $\Uglambda$ instead of $\pr_{\ker(\chi_\lambda)}$ and $\Ugmod_{\ker(\chi_\lambda)}$, respectively. Given $\lambda,\mu \in \frt^*_E$ such that $\nu := \mu -\lambda$ is integral, the translation functor $\Tlamu$ is the defined by

\begin{equation}\label{intro T functor}
\begin{array}{lccc}\Tlamu: & \Ugzfin & \ra & \Ugzfin\\
&&&\\
& M & \rightsquigarrow & \pr_{|\mu|} (L(\onu) \ot_E \pr_{|\lambda|}(M)) \;,
\end{array}\tag{1}
\end{equation}

\vskip8pt

where $\onu$ is the dominant weight in the (linear) Weyl orbit of $\nu$ and $L(\onu)$ the irreducible finite-dimensional module with highest weight $\onu$. Under the assumptions (i)-(iii) in Thm. 1 below, the functor $\Tlamu$ induces an equivalence of categories  $\Uglambda \ra \Ugmu$.\footnote{This also holds when condition (i) is replaced by the weaker condition that $\mu-\lambda$ is integral.}

\vskip8pt

{\it Translation functors for locally analytic representations.} Analogous translation functors can be introduced in other areas of representation theory, and our intention is to carry this out in the framework of locally analytic representations of $p$-adic groups, as introduced and studied by P. Schneider and J. Teitelbaum \cite{S-T1,S-T2}. 

\vskip8pt

Here we do not strife for utmost generality and consider a split reductive group $\bG$ over $F$ (itself a finite extension of $\Qp$), a Borel subgroup $\bB$ with Lie algebra $\frb$ and a maximal split torus $\bT \sub \bB$ with Lie algebra $\frt$. The corresponding groups of $F$-valued points will be denoted by $G$, $B$, $T$. Instead of working directly with locally $F$-analytic representations of $G$, it is more convenient to work with modules over the locally analytic distribution algebra $\DG$ of $G$, which is the continuous dual space of the space of locally $F$-analytic $E$-valued functions on $G$. There is a canonical injective algebra homomorphism $\UgE \hra \DG$. Important for our purposes is that $\frz_E$ is mapped into the center of $\DG$ via this map, so that we have the notion of $\frz$-finite $\DG$-modules available. The translation functor 

\begin{equation}\label{intro cT functor}
\begin{array}{lccc}\cTlamu: & \Dzfin & \ra & \Dzfin\\
&&&\\
& M & \rightsquigarrow & \pr_{|\mu|} (L(\onu) \ot_E \pr_{|\lambda|}(M))
\end{array}\tag{2}
\end{equation}

\vskip8pt

is then defined exactly as in (\ref{intro T functor}), except that we require that $\onu$ lifts to an algebraic character of $\bT$, which in turn ensures that $L(\onu)$ lifts to an algebraic representation of $\bG$. While those functors on $\DG$-modules are defined by the same formula as for $\UgE$-modules, we use a different font to distinguish them. The following theorem concerns the category of coadmissible $\DG$-modules as defined in \cite{S-T2}.  

\vskip8pt

{\bf Theorem 1.}\label{main1} (Thm. \ref{D-equivalence}) {\it Let $\lambda, \mu \in \frt^*_E$ satisfy the following conditions:

\begin{itemize}
    \item[(i)] $\mu - \lambda$ lifts to an algebraic character of $\bT$.
    \item[(ii)] $\lambda$ and $\mu$ are both anti-dominant, cf. \ref{anti-dominant}.
    \item[(iii)] The stabilizers in $W$ of $\lambda$ and $\mu$ for the dot-action are the same.
\end{itemize}

\vskip8pt

Then $\cTlamu$ induces an equivalence of categories $\Dlambda \ra \Dmu$ with quasi-inverse $\cTmula$.} \qed

\vskip8pt

The proof of this result is by reduction to the categorical equivalences of \cite{BeGe} for $\UgE$-modules. 

\vskip8pt

{\it Relations with category $\cO$.} Translation functors on $\Ugzfin$ preserve the Bernstein-Gelfand-Gelfand category $\cO$ and furnish important information about this category. For a standard-parabolic subgroup $\bP$ with Lie algebra $\frp$ denote by $\cO^\frp$ the category of modules $M$ in $\cO$ on which $\frp$ acts locally finite. Following  \cite{O-S2,O-S3} and \cite{AS} we consider here the category $\cO^{P,\infty}$ of pairs $(M,\phi)$, where $M$ is an object of the category $\cO^{\frp,\infty}$, which is the extension-closure of $\cO^\frp$, and where $\phi: P=\bP(F) \ra \End_E(M)^\x$ is a lift of the $\frp$-action to a locally $F$-analytic action of $P$. Generalizing the work done in these papers there is an exact functor

\[\chcFGP: \cOPinfty \times \Rep^{\sm}_E(L_P)^\sadm \lra  \DGcoadzfin \;,\]

\vskip8pt

where $\Rep^{\sm}_E(L_P)^\sadm$ is the category of strongly-admissible smooth representations of the Levi subgroup $L_P$ of $P$. This functor commutes with the translation functors in (\ref{intro T functor}), extended to $\cO^{P,\infty}$ in a straightforward way, and (\ref{intro cT functor}) in the following sense.

\vskip8pt

{\bf Theorem 2.} (Thm. \ref{F-T commute}) {\it Let $\lambda, \mu \in \frt^*_E$ be such that $\mu -\lambda$ is algebraic. For any $(M,V)$ in $\cOPinfty \times \Rep^{\sm}_E(L_P)^\sadm$ there is a canonical isomorphism

\[ \cTlamu\Big(\chcFGP(M,V)\Big) \cong \chcFGP \Big(T_\lambda^\mu(M),V\Big)\]

\vskip8pt

which is natural in $M$ and $V$.}\qed

\vskip8pt

As a corollary of this result we deduce some formulas for the effect of the translation functors on locally analytic principal series representations (or rather their coadmissible modules), using how translation functors act on Verma modules in category $\cO$ \cite[sec. 7]{Hu}.

\vskip8pt

{\it Relations with the p-adic Langlands correspondence for $\GL_2(\Qp)$.} As an application of the methods initiated in this paper, we consider the case of $\GL_{2,\Qp}$ and how translation functors act on locally analytic representations of the form $\Pi(V)^\la$ (or rather their coadmissible modules). Here, $V$ is a 2-dimensional absolutely irreducible representation of the absolute Galois group of $\Qp$ and $\Pi(V)^\la$ denotes the associated locally analytic representation\footnote{The associated locally analytic representation consists of the subspace of locally analytic vectors, but this subspace is equipped with an intrinsic topology which makes it complete and is finer than the subspace topology of $\Pi(V)$.} of the Banach space representation $\Pi(V)$ of the $p$-adic local Langlands correspondence \cite{Colmez_GL2,ColDoPa}.  

\vskip8pt

Our result in this direction concerns absolutely irreducible trianguline Galois representations $V(s)$, where the parameter $s = (\delta_1, \delta_2, \infty)$ in the trianguline variety $\sS_{\irr}$ is generic in the sense that $(\delta\delta_2^{-1})(x) \neq x^i|x|$ for any $i \in \bbZ_{\ge 1}$. In the following we write $\Pi(s)$ for the Banach space representation attached to the trianguline Galois representation $V(s)$. 

\vskip8pt

{\bf Theorem 3.} (Thm. \ref{trianguline thm}) {\it Let $s$ be a generic point of $\sS_{\irr}$. Let $\wnu$ be an algebraic character of $\bT$ and $\theta$ a locally analytic character of $\Qpx$ satisfying \ref{cond theta}. Define $\wlambda = \wlambda(s)$ as in \ref{lambdatilde} and $\wmu$, $\whs$ as in \ref{hat defs}. Assume \ref{cond1}, \ref{cond3}, and \ref{cond2} are fulfilled. Then:}

\vskip8pt

\begin{enumerate}
\item[(i)] $\whs$ {\it is a generic point in $\sS_{\irr}$.}

\vskip5pt

\item[(ii)] {\it There is an isomorphism of $\DG$-modules}

\[(\theta \ot \cTlamu)\Big((\Pi(s)^{\la})'\Big) \cong \Big(\Pi(\whs)^{\la}\Big)' \;,\] 

\vskip8pt

where $\lambda = \rmd \wlambda$ and $\mu = \rmd \wmu$. \qed
\end{enumerate}

\vskip25pt

While this result seems to be very technical and apparently only holds under strong assumptions, we would like to point out that those assumptions are quite natural and allow for many examples where translation functors move representations of the form $\Pi(s)^{\la}$ to representations of the same kind. 

\vskip8pt

{\it Acknowledgments.} 
We thank Andrés Sarrazola Alzate for carefully reading the manuscript and giving helpful suggestions. A.L. would like to thank the staff of Indiana University Mathematics department for all their help during an eventful year.
\vskip8pt

{\it Notation.} We let $F/\Qp$ be a finite field extension of $\Qp$, and we let $E$, the `coefficient field', be a finite field extension of $F$. 

\vskip8pt

Let $\bG$ be a connected split reductive algebraic group over $F$, and $\bB \supset \bT$  a Borel subgroup and a maximal split torus, respectively. We denote by $G \supset B \supset T$ the groups of $F$-valued points of $\bG$, $\bB$, and $\bT$, respectively. $H \sub G$ will always denote a compact open subgroup of $G$. We let $\frg \supset \frb \supset \frt$ be the Lie algebras of $G$, $B$ and $T$, respectively. Denote by $\Phi = \Phi(\frg,\frt)$ the associated root system, and by $\Phi^+ \sub \Phi$ the set of positive roots with respect to $\frb$. 

\vskip8pt

Let $\frz  = Z(U(\frg))$ be the center of the universal enveloping algebra of $\frg$. The subscript `$E$' denotes always the base change to $E$. For example $\frz_E = \frz \ot_F E$ is the center of $\UgE$. The Weyl group $W$ of the pair $(\frg,\frt)$ is, by definition, the Weyl group of the pair $(\frg',\frt \cap \frg')$, where $\frg' = [\frg,\frg]$ is its derived subalgebra $[\frg,\frg]$. The group $W$ acts trivially on the dual space of the center $Z(\frg)$ of $\frg$. The `dot-action' of $W$ on weight space $\frt^*_E  = \Hom_F(\frt,E)$ is defined by $w \cdot \lambda := w(\lambda +\rho) -\rho$, where $\rho = \frac{1}{2}\sum_{\alpha \in \Phi^+} \alpha$. For the linear $W$-action on weight space we use parentheses, e.g., $w(\nu)$ denotes the $E$-linear action of $w \in W$ on $\nu \in \frt^*_E$. By $|\lambda|$ we denote the orbit of $\lambda$ under the dot-action of $W$. 

\vskip8pt

Given $\lambda \in \frt^*_E$ we denote $E_\lambda$ the 1-dimensional $\frt_E$-module $E$ with the property that $h.1 = \lambda(h)$ for all $h \in \frt_E$. We extend $E_\lambda$ to a $\frb_E$-module via the map $\frb_E \ra \frb_E/[\frb,\frb]_E = \frt_E$. Then $\chi_\lambda: \frz \ra E$ denotes the character by which $\frz_E$ acts on the Verma module $M(\lambda) = \UgE \ot_{\UbE} E_\lambda$ with highest weight $\lambda$. One has $\chi_\lambda = \chi_\mu$ if and only if $|\lambda| = |\mu|$, cf. \cite[4.115]{KnappVogan_Cohomological}.\footnote{Note that in \cite{KnappVogan_Cohomological} the character $\chi_\lambda$ is that of the Verma module with highest weight $\lambda - \rho$, and that the Weyl group orbits are with respect to the linear action.}

\vskip8pt

$\lambda \in \frt^*_E$ is called {\it integral} if $\pair{\lambda}{\alpha^\vee} \in \Z$ for all $\alpha \in \Phi$. By $\Lambda$ we denote the set of integral weights, by $\Lambda^+ = \{\lambda \in \frt^*_E \midc \forall \alpha \in \Phi^+: \pair{\lambda}{\alpha^\vee} \in \bbZ_{\ge 0}\}$ the set of dominant integral weights, and by $\Lambda_r = \langle \Phi \rangle_\Z \sub \Lambda$ the lattice generated by the roots.

\vskip8pt

If $X$ is a topological vector space over $E$, we denote by $X'_b$ the space of continuous linear forms on $X$, equipped with the strong topology \cite[p. 35]{Schneider_NFA}. Only for $\frt$ or $\frt_E$ we write $\frt^*$ and $\frt^*_E$, respectively, for the dual space (in order to conform to standard conventions).  

\vskip8pt

All modules are considered to be left modules, and $\Rmod$ denotes the category of left $R$-modules.

\vskip8pt

\section{Translation functors}

\subsection{Locally analytic distribution algebras and related rings}\label{prelim}
We recall some concepts and results from the paper \cite{S-T1,S-T2} by P. Schneider and J. Teitelbaum. 

\vskip8pt

\begin{para}{\it Locally analytic distribution algebras.}
Let $\DG := C^\la(G,E)'_b$ be the locally analytic distribution algebra, which is the strong dual of the space of $E$-valued locally $F$-analytic functions on $G$.\footnote{Because the coefficient field $E$ will remain the same throughout the entire paper, we do not indicate the dependence of $E$ in $\DG$ and hence do not use the notation $D(G,E)$ as in \cite{S-T1,S-T2}.} The multiplication on $\DG$ is given by 
\[(\delta_1 \delta_2)(f) = \delta_1\Big(g \mapsto \delta_2(h \mapsto f(gh))\Big) \;,\]
for $\delta_1,\delta_2 \in \DG$ and $f \in C^\la(G,E)$. Given a locally $F$-analytic representation $V$ of $G$ on an $E$-vector space of compact type, the strong dual space $V'_b$ is naturally a module over $\DG$: for $\delta \in \DG$ and $\mu \in V'_b$ one defines $(\delta.\mu)(v) = \delta(g \mapsto \mu(g^{-1}.v))$.  By \cite[3.3]{S-T1} the functor $V \ra V'_b$ is an anti-equivalence between the category of locally analytic $G$-representations on $E$-vector spaces of compact type with continuous linear $G$-maps to the category of separately continuous $\DG$-modules on nuclear Fr\'echet spaces with continuous $\DG$-module maps. 
\end{para}

\vskip8pt

\begin{para}{\it Embedding of the universal enveloping algebra and of $E[G]$.}
Any $\frx \in \frg$ gives rise to a locally analytic distribution by 
\[f \mapsto \frac{d}{dt}\Big(t \mapsto f(e^{-t\frx})\Big)|_{t=0}\;,\] 
which in turn induces an injective morphism of $E$-algebras $\UgE \hra \DG$ \cite[p. 449]{S-T1}. Under this map the center $\frz_E$ of $\UgE$ is mapped into $Z(\DG)$, where $Z(\DG)$ denotes the center of the distribution algebra of $G$ \cite[3.7]{S-T1}.

\vskip8pt

For $g \in G$ we denote by $\delta_g \in \DG$ the delta-distribution at $g$ which is defined by $\delta_g(f) = f(g)$. The resulting map $G \ra \DG^\times$ is a group homomorphism and extends $E$-linearly to a morphism of $E$-algebras $E[G] \ra \DG$ which has dense image \cite[3.1]{S-T1}. 
\end{para}

\vskip8pt

\begin{para}{\it Schneider-Teitelbaum completions and coadmissible modules.}
For any compact open subgroup $H$ of $G$ and $r \in (\frac{1}{p},1)$ there exist completions $\DrH = D_r(H,E)$ of $\DH$ with respect to certain algebra semi-norms $q_r$ such that 
\[\DH = \varprojlim_{\frac{1}{p}< r < 1} \DrH \;.\] 
The rings $\DrH$ are noetherian Banach algebras and this gives $\DH$ the structure of a Fr\'echet-Stein algebra \cite[sec. 4 and Thm. 5.1]{S-T2}. As for any Fr\'echet-Stein algebra, we can consider the category of coadmissible $\DH$-modules \cite[p. 152]{S-T2}. A coadmissible $\DH$-module has a canonical topology \cite[p. 155]{S-T2}. A $\DG$-module $M$ is called coadmissible, if $M$ is coadmissible as a $\DH$-module for some (equivalently, any) compact open subgroup $H$ of $G$. The canonical topology on a coadmissible $\DG$-module $M$, considered as a $\DH$-module, does not depend on the compact open subgroup $H$. It is a nuclear Fr\'echet space and gives $M$ the structure of a separately continuous $\DG$-module \cite[p. 174]{S-T2}. This topology is called the canonical topology on $M$. Any morphism between coadmissible $\DG$-modules is continuous for the canonical topologies. We denote the category of coadmissible $\DG$-modules by $\DGcoad$. 
\end{para}

\vskip8pt

\begin{para}{\it The ring $\DgP$.}
Let $\bP \sub \bG$ be a standard parabolic subgroup, and set $P = \bP(F)$ and $\frp = \Lie(\bP)$. We denote by $\DgP \sub \DG$ be the smallest subring containing $\UpE$ and $\DP$. It follows from \cite[3.5]{O-S2} that as a $\UgE\mbox{-}\DP$-bimodule the multiplication map induces an isomorphism 

\begin{numequation}\label{U-DP-bimodule iso}
\UgE \ot_{\UpE} \DP \stackrel{\simeq}{\lra} \DgP  \;,
\end{numequation}

cf. also \cite[4.1]{ScSt}.
\end{para}

\vskip8pt

\subsection{Categories of \texorpdfstring{$\frz$}{}-finite modules}\label{zfinite} In the following we let $R$ be a unital $\frz_E$-algebra, and we denote by $\Rmod$ the category of left $R$-modules.  

\begin{para}{\it Categories of $\frz$-finite $R$-modules.}\label{Rzfin modules} 
Given $M$ in $\Rmod$ and an ideal $J \sub \frz_E$, we set $M_{J = 0} = \{m \in M \midc J.m = \{0\}\}$. This is an $R$-submodule of $M$.  
By $\Rzfin$ we denote the full subcategory of $\Rmod$ consisting of objects $M$ such that $M = \bigcup_J M_{J=0}$, where $J$ runs through the ideals of finite codimension in $\frz_E$. We call those modules $\frz$-finite. Given a maximal ideal $\frm \sub \frz_E$, we let $\Rfrm$ be the full subcategory of $\Rzfin$ consisting of modules $M$ with $M = \bigcup_{n \ge 1} M_{\frm^n=0}$.

\vskip8pt

If $J \subset \frz_E$ is of finite codimension, then $\frz_E/J$ is an Artin ring. By the general theory of Artin rings, $J = 
\bigcap_{i=1}^r J_i$ with primary ideals $J_i$ whose radicals $
\sqrt{J_i}$ are maximal and pairwise different \cite[sec. 8]{AtiyahMacDonald}. Then $\frz_E/J = \prod_{i=1}^r \frz_E/J_i$. If $M = M_{J=0}$ we have 
\[M = \frz_E/J \otimes_{\frz_E} M = \prod_{i=1}^r \frz_E/J_i \otimes_{\frz_E} M \;.\]
Notice that if $\frm = \sqrt{J_i}$, then $\frz_E/J_i = (\frz_E/J)_\frm = (\frz_E)_\frm \otimes \frz_E/J$, and thus \[\frz_E/J_i \otimes_{\frz_E} M = (\frz_E)_\frm \otimes_{\frz_E} M = M_\frm\;,\] 
where $M_\frm$ is an object of $\Rfrm$. Given $M$ in $\Rzfin$ we thus have 

\[M  = \; \bigcup_J M_{J=0} 
\; = \; \bigcup_J \bigoplus_{\frm \in \Max(\frz_E)} (M_{J=0})_\frm = \bigoplus_{\frm \in \Max(\frz_E)} M_\frm \;,\]

\vskip8pt

where $J$ runs through the ideals of finite codimension in $\frz_E$. Denote by $\pr_\frm$ the canonical map $M \ra M_\frm$; it is an exact endo-functor on $\Rzfin$. Occasionally, we consider $\pr_\frm$ as a functor from $\Rzfin$ to $\Rfrm$. It should be clear from the context when $\pr_\frm$ is considered in the latter sense.  
\end{para}

\vskip8pt

The following assertions are straightforward. 

\vskip8pt

\begin{lemma}\label{Homzeroes}
Let $M,N\in \Rzfin$ and $\frm, \frn \in \Max(\frz_E)$.
\vskip8pt
\begin{enumerate}
\item[(i)] $\pr_\frm$ is an exact endo-functor of $\Rzfin$. \vskip5pt
\item[(ii)] If $\frm \neq \frn$, then $\Hom_R(M_\frm, N_\frn)= 0$. \vskip5pt
\item[(iii)] $\Hom_R(M, N_{|\lambda|}) = \Hom_R(M_{|\lambda|}, N)= \Hom_R(M_{|\lambda|}, N_{|\lambda|})$. \vskip5pt
\item[(iv)] If $R \ra S$ is a morphism of $\frz_E$-algebras, then $\pr_\frm(S \ot_R M) = S \ot_R \pr_\frm(M)$.
\end{enumerate}
\end{lemma}

\vskip8pt

Recall that for $\lambda\in \frt^*_E$ we denote by $\chi_\lambda: \frz_E \ra E$ the character by which $\frz_E$ acts on the Verma module $M(\lambda)$ with highest weight $\lambda$. Put $\frm_\lambda = \ker (\chi_\lambda)$, which only depends on $|\lambda|$. If $\frm \in \Max(\frz_E)$ is equal to $\frm_\lambda$, then we write $M_{|\lambda|}$, $\Rlambda$, and $\pr_{|\lambda|}$, instead of $M_{\frm_\lambda}$, $\Rmod_{\frm_\lambda}$, and $\pr_{\frm_\lambda}$. 

\vskip8pt

\begin{para}{\it Categories of $\frz$-finite modules over distribution algebras.}\label{Dzfin modules} 
We apply the previous notation and remarks to the case when $R$ is $\UgE$, $\DG$ od $\DrH$. As mentioned in \ref{prelim}, there is a canonical morphism $\UgE \ra \DG$ of $E$-algebras, which has the property that the image of $\frz_E \sub \UgE$ is mapped into the center of $\DG$, and the same applies to $\DrH$ as well.

\vskip8pt

Finally, the superscript “fg” indicates the full subcategory of finitely generated modules in the corresponding category, e.g., $\DrHfglambda$ consists of finitely generated modules in $\DrHlambda$.
\end{para}

\vskip8pt

\subsection{Tensoring with finite-dimensional representations}\label{tensoring}

\begin{para}{\it A preliminary remark.} The diagonal map $G \ra G \times G$ induces a continuous morphism of $E$-algebras 
\[\DG \ra D(G \times G)\] 
and the latter ring is canonically isomorphic to the injective complete tensor product $\DG \wot_{E,\iota} \DG$ \cite[p. 312]{ST_duality}. This shows that $\DG$ does not have a co-multiplication in the algebraic sense (i.e., a ring homomorphism to $\DG \ot_E \DG$), unless $G$ is discrete (in which case $\DG = E[G]$). In particular, a tensor product $L \ot_E M$ of $\DG$-modules $L$ and $M$ does not always carry the structure of a $\DG$-modules. However, if one of the modules is a finite-dimensional locally $F$-analytic representation, then this is not a problem, as is stated below. We recall that any locally $F$-analytic representation of $G$ carries the structure of a separately continuous $\DG$-module \cite[3.2]{S-T1}.
\end{para}

\begin{lemma}\label{A.2.3} Let $R$ be any of the rings $\DG$, $\DrH$, or $\DgP$. Suppose $L$ is a finite-dimensional locally $F$-analytic representation of $G$. Then, for any $R$-module $X$, there exists a natural $R$-module structure on the tensor product $L\ot_E X$ with the property that for all $g \in G$ (resp. $g \in H$, resp. $g \in P$), $l \in L$, and $x \in X$ one has

\[\delta_g.( l\ot x) = (g.l) \ot (\delta_g.x) \;.\]
\end{lemma}

\vskip8pt

\begin{proof}
For $R = \DG$ this is \cite[6.2.1]{AS}. By restricting to $\DgP \sub \DG$, this implies the assertion for $R = \DgP$. A quick inspection shows that the proof given in \cite[6.2.1]{AS} stays valid when $\DrH$-modules are being considered.
\end{proof}

\vskip8pt

\begin{lemma}\label{finite adjunction}
Let $L$ be a finite-dimensional locally $F$-analytic representation of $G$ over $E$. Let $R$ be one of the rings $\UgE$, $\DG$, $\DrH$, or $\DgP$. 

\vskip8pt

For any $R$-modules $M$ and $X$ we have a natural isomorphism of $E$-vector spaces

\[\Hom_R (L\ot_E M, X) \cong \Hom_R(M, \Hom_E(L, X)) \;.\]
\end{lemma}

\vskip8pt

\begin{proof} For $R = \UgE$ this statement is \cite[Cor 4.3]{Kn}. For $R = \DG$, this has been shown in \cite[6.3.1]{AS}, and this implies the assertion for $R = \DgP$. Again, a quick inspection shows that the proof of \cite[6.3.1]{AS} stays valid when $\DrH$-modules are being considered. 
\end{proof}

\vskip8pt

\begin{theorem}\label{finite commuting}
Let $L$ be a finite-dimensional locally $F$-analytic representation of $G$ over $E$. Let $S$ be of the rings in \ref{finite adjunction}. Let $R$ be of the rings in \ref{finite adjunction} which is contained in $S$, or $R = \DH$ if $S = \DrH$. Let $M$ be an $R$-module. Then
there exists an isomorphism of $S$-modules which is natural in $L$ and $M$

\[S \ot_R (L\ot_E M)  \lra  L \ot_E (S \ot_R M) \;.\]
\end{theorem}

\vskip8pt

\begin{proof}
Via the Yoneda lemma, the assertion is equivalent to showing that for any $S$-module $X$ there is a natural isomorphism

\[\begin{array}{rcl}
\Hom_S(S \ot_R (L\ot_E M), X) & \cong & \Hom_R(L\ot_E M, X)\\
&&\\
& \stackrel{\ref{finite adjunction}}{\cong} & \Hom_R( M, \Hom_E(L, X)) \\
&&\\    
& \cong & \Hom_{S}(S \ot_R M, \Hom_E(L, X))\\
&&\\    
&\stackrel{\ref{finite adjunction}}{\cong} & \Hom_S\left(L \ot_E (S \ot_R M), X\right) \;. 
\end{array}\]
\end{proof}

\vskip8pt

\begin{lemma}\label{zfinite stable} Let $\cA$ be one of the categories $\Dzfin$, $\DGcoad$,\linebreak $\DGcoadzfin$, $\DrHzfin$, $\DrHmod^{\rm fg}$, or $\DrHfgzfin$. 

\vskip8pt

For any $M$ belonging to $\cA$ and any finite-dimensional locally $F$-analytic representation $L$ of $G$, the tensor product $L \ot_E M$ belongs to $\cA$ as well.  
\end{lemma}

\begin{proof} (i)  For $\frz$-finiteness we only need to consider the underlying structure as a $\UgE$-module. \cite[2.6 (ii)]{BeGe} thus implies that $L \ot_E M$ is $\frz$-finite.

\vskip8pt

(ii) Next we show that $\DrHmod^{\rm fg}$ is stable under tensoring with $L$. Given $M$ in $\DrHmod^{\rm fg}$, let $N \sub M$ be a finitely generated $\UgE$-submodule which generates $M$ as $\DrH$-module. Then $L \ot_E N$ is finitely generated as $\UgE$-module \cite[2.3 (i)]{BeGe}, and by \ref{finite commuting} the module $L\ot_E(\DrH \ot_{\UgE} N)$ is finitely generated as $\DrH$-module. Since $L \ot_E M$ is a quotient of the latter module, it is finitely generated as well.

\vskip8pt 

(iii) Now let $M$ be a coadmissible $\DG$-module. Then $M_r := \DrH\ot_\DH M$ is finitely generated, and so is $L\ot_E(\DrH\ot_\DH M)$, by what we have just seen. Using \ref{finite commuting}, we see that $(L \ot_E M)_r := \DrH\ot_\DH (L\ot_E M) \cong L \ot_E M_r$ is finitely generated. Because $M$ is the projective limit of the $M_r$, the projective limit of the $(L \ot_E M)_r$ is $L \ot_E M$, and $L \ot_E M$ is thus coadmissible as $\DH$-module. 
\end{proof}

\vskip8pt

\subsection{Translation functors: definition and basic properties} 

\begin{para}\label{Tlamu}{\it Translation functors for $\UgE$-modules.} We begin by recalling the translation functors as considered in \cite{BeGe} and \cite[7.6]{Hu}, following the setup and notation of the latter reference.\footnote{Note that in \cite[7.6]{Hu} these functors are only introduced in the context of category $\cO$. But their definition makes sense more generally for $\frz$-finite $\UgE$-modules.} 

\vskip8pt

Let $\lambda, \mu \in \frt^*_E$ be such that $\nu := \mu-\lambda$ is integral, i.e., $\pair{\nu}{\alpha^\vee} \in \Z$ for all $\alpha \in \Phi$. Let $\onu \in \Lambda^+$ be the unique dominant weight in the $W$-orbit of $\nu$ for the linear Weyl group action, and let $L(\onu)$ be the corresponding finite-dimensional absolutely irreducible $\frg_E$-representation. Then the  
translation functor

\[\Tlamu: \Ugzfin \lra \Ugzfin\] 
is defined by
\[\Tlamu(M)=\pr_{|\mu|}\Big(L(\onu)\ot_E \pr_{|\lambda|}(M)\Big) \;.\]

\vskip8pt

This functor is oftentimes considered as a functor $\Uglambda \ra \Ugmu$. It will be clear from the context when this is the case.  
\end{para}

\vskip8pt

\begin{para}{\it Algebraic weights.}
$\lambda \in \frt^*_E$ is said to be {\it algebraic} (for $\bT$) if it is in the image of the (injective) map 

\[X^*(\bT) \ra \frt^*_E \;, \;\; \tau \mapsto {\rm d}\tau \;.\]

\vskip8pt

Algebraic weights are always integral, but the converse is not necessarily the case, even if $\bG$ is semisimple.\footnote{If $\bG$ is semisimple and simply-connected, however, then the algebraic and integral weights are the same \cite[II, 1.6]{JantzenRepresentations}.} Weights in $\Lambda_r$ are always algebraic. 
\end{para}

\begin{para}\label{compatible}{\it Compatible weights.} We call $\lambda, \mu \in \frt^*_E$ {\it compatible} if $\lambda -\mu$ is algebraic. In this case, we let $\onu$ be the unique dominant weight in the $W$-orbit of $\nu$ with respect to the linear Weyl group action. The irreducible finite-dimensional representation $L(\onu)$ with highest weight $\onu$ lifts then (uniquely) to an algebraic representation of $\bG$ \cite[II, 1.20]{JantzenRepresentations}. We consider $L(\onu)$ as a locally analytic representation of $G$ and hence as a $\DG$-module \cite[3.2]{S-T1}. 
\end{para}

\begin{rem} Because of the possible discrepancy of integral and algebraic weights, the concept of compatible weights is not the same as that introduced in \cite[ch. 7]{Hu}. Compatible weights, as introduced here, are compatible in the sense of \cite[ch. 7]{Hu}, but the converse does not necessarily hold.
\qed\end{rem}

\begin{dfn} Let $\lambda, \mu \in \frt^*_E$ be compatible in the sense of \ref{compatible}. We define the translation functor

\[\cTlamu:\Dzfin \ra \Dzfin\] 
by
\[\cTlamu(M)=\pr_{|\mu|}\Big(L(\onu)\ot_E \pr_{|\lambda|}(M)\Big) \;.\]
\end{dfn}

\vskip8pt

Note that this is well-defined by \ref{zfinite stable}. Furthermore, again by \ref{zfinite stable}, the functor $\cTlamu$ induces functors $\cTclamu$ and $\cTrlamu$ on the categories $\DGcoadzfin$ and $\DrHfgzfin$, respectively. Oftentimes we consider those functors as functors 

\[\begin{array}{lcl}
\Dlambda & \lra &\Dmu \;,\\ 
&&\\
\DGcoadlambda & \lra &\DGcoadmu \;,\\
&&\\
\DrHfglambda &\lra &\DrHfgmu \;.
\end{array}\]

\vskip8pt

\begin{lemma}\label{translation exact}
Let $\lambda, \mu \in \frt^*_E$ be compatible. Then $\cTlamu$ is an exact functor. 
\end{lemma}

\begin{proof}
Taking the tensor product over a field is exact. The projection operators $\pr_{|\lambda|}$ and $\pr_{|\mu|}$ are exact by \Cref{Homzeroes}. $\cTlamu$ is a composition of these functors and hence exact.
\end{proof}

\vskip8pt

\begin{theorem}\label{T adjoint} The functor $\cTlamu$ is left and right adjoint to $\cTmula$. That is, for all $M, N\in \Dzfin$ there is an isomorphism 

\[\Hom_\DG(\cTlamu M, N)\cong \Hom_\DG(M, \cTmula N) \;,\]

\vskip8pt

which is natural in $M$ and $N$. The same holds when $\lambda$ and $\mu$ are exchanged.

\vskip8pt

Similarly, $\Tlamu$ ($\cTrlamu$, $\cTclamu$, resp.) is left and right adjoint to $\Tmula$ ($\cTrmula$, $\cTcmula$, resp.) on $\Ugzfin$ ($\DrHzfin$, $\DGcoadzfin$, resp.).
\end{theorem}

\vskip8pt

\begin{proof}
Let $w_\circ$ be the longest Weyl group element. Denote by $L(\onu)'$ the dual space of $L(\onu)$, equipped with the contragredient group action. Then $L(\overline{\nu})' \cong L(-w_\circ(\onu))$ \cite[1.6]{Hu}. Moreover, $-w_\circ(\onu)$ is the unique $W$-conjugate in $\Lambda^+$ of $\lambda-\mu$ \cite[7.2]{Hu}.  Let $R$ be one of the rings $\UgE$, $\DG$, or $\DrH$. Then:

\vskip8pt

\[\begin{array}{rcl}
\Hom_R(\cTlamu M, N) & \stackrel{\ref{Homzeroes}}{\cong} & \Hom_R(\cTlamu M, \pr_{|\mu|}(N))\\
&&\\
& \stackrel{\ref{Homzeroes}}{\cong} & \Hom_R(L(\onu)\ot_E \pr_{|\lambda|}(M), \pr_{|\mu|}(N))\\
&&\\
& \stackrel{\ref{finite adjunction}}{\cong} & \Hom_R(\pr_{|\lambda|}(M), L(\onu)^*\ot_E \pr_{|\mu|}(N))\\
&&\\
& \cong & \Hom_R(\pr_{|\lambda|}(M), L(-w_\circ\onu)\ot_E \pr_{|\mu|}(N))\\
&&\\
& \stackrel{\ref{Homzeroes}}{\cong} & \Hom_R(M, \pr_{|\lambda|}(L(-w_\circ\onu)\ot_E \pr_{|\mu|}(N)))\\ 
&&\\
& = & \Hom_R(M, \cTmula N) \;.
\end{array}\]
\end{proof}   

\vskip8pt

\begin{lemma}\label{T commute basechange}
Let $\lambda, \mu \in \frt^*_E$ be compatible. Let $S$ be one of the rings in \ref{finite adjunction}. Let $R$ be one of the rings in \ref{finite adjunction} which is contained in $S$, or $R = \DH$ if $S = \DrH$. Denote by $\cT^\mu_{R,\lambda}$ and $\cT^\mu_{S,\lambda}$ the translation functors on $R\mbox{-}{\rm mod}_{\frz\mbox{-}{\rm fin}}$ and $S\mbox{-}{\rm mod}_{\frz\mbox{-}{\rm fin}}$, respectively. Then there is a canonical isomorphism
\[S \otimes_R \cT^\mu_{R,\lambda}(-)  \lra \cT^\mu_{S,\lambda} \Big(S \otimes_R (-)\Big)\]
of functors
\[R\mbox{-}{\rm mod}_{\frz\mbox{-}{\rm fin}} \lra S\mbox{-}{\rm mod}_{\frz\mbox{-}{\rm fin}} \;.\]
\end{lemma}

\begin{proof}
\[\begin{array}{rcl}
S \otimes_R \cT^\mu_{R,\lambda}(M) & = &  S \otimes_R \pr_{|\mu|}(L(\onu) \ot_E \pr_{|\lambda|}(M)\\
&&\\     
& \stackrel{\ref{Homzeroes}}{\cong} & \pr_{|\mu|}\Big(S \ot_R (L(\onu) \ot_E \pr_{|\lambda|}(M))\Big)\\
&&\\
& \stackrel{\ref{finite commuting}}{\cong} & \pr_{|\mu|}\Big(L(\onu) \ot_E (S \ot_R \pr_{|\lambda|}(M))\Big)\\
&&\\
& \stackrel{\ref{Homzeroes}}{\cong} & \pr_{|\mu|}\Big(L(\onu) \ot_E \pr_{|\lambda|}(S \ot_R M)\Big)\\
&&\\
& = & \cT^\mu_{S,\lambda}\Big(S \otimes_R M\Big) \;.
\end{array}\]
\end{proof}

\section{Categorical Equivalences}
The main result of this section is that the functors $\cT^\mu_{c,\lambda}$ are equivalences of categories, under suitable assumptions on $\lambda$ and $\mu$. Our strategy to prove this result is by reduction to the case of $\UgE$-modules. We begin by recalling the setup and results from \cite[ch. 7]{Hu} and \cite{BeGe}. With regard to notation and conventions (e.g., dot-action of the Weyl group versus linear action, and twisted Harish-Chandra homomorphism versus untwisted Harish-Chandra homomorphism) we follow \cite{Hu}.

\vskip8pt

\subsection{Equivalences for categories of \texorpdfstring{$\Ug$-modules}{}}
Recall that $\lambda \in \frt^*_E$ is called {\it anti-dominant}\footnote{as defined in \cite[3.5]{Hu}} if

\begin{numequation}\label{anti-dominant}
\mbox{for all } \alpha \in \Phi^+: \pair{\lambda+\rho}{\alpha^\vee} \notin \Z_{\ge 1} \;.
\end{numequation}

For $\lambda\in \frt^*_E$ let $W^\circ_\lambda=\{w\in W\st w \cdot \lambda = \lambda\}$ be the stabilizer of $\lambda$ for the dot-action.\footnote{This is the notation used in \cite[7.4]{Hu} and \cite{JantzenModuln}.}  

\vskip8pt

\begin{theorem}\label{BeGe theorem} Let $\lambda, \mu \in \frt^*_E$ satisfy the following conditions:

\begin{itemize}
    \item[(i)] $\lambda$ and $\mu$ are compatible in the sense of \ref{compatible}.\footnote{The theorem holds under the weaker assumption that $\mu - \lambda$ is integral. As we will later on need the stronger assumption that they are compatible, we assume this here too.}
    \item[(ii)] $\lambda$ and $\mu$ are both anti-dominant.
    \item[(iii)] $W^\circ_\lambda = W^\circ_\mu$.
\end{itemize}

\vskip8pt

Then the functors $\Tlamu$ and $\Tmula$ from \ref{Tlamu} induce equivalences of categories

\[\Tlamu: \Uglambda \lra \Ugmu\] 
and 
\[\Tmula: \Ugmu \lra \Uglambda \;.\]

\vskip8pt

They also induce equivalences of categories $\Ugfglambda \ra \Ugfgmu$ and $\Ugfgmu \ra \Ugfglambda$. One has natural isomorphisms $\Tlamu \circ \Tmula \cong \id$ and $\Tmula \circ \Tlamu \cong \id$.
\end{theorem}

\begin{proof} These assertions will be deduced from \cite[4.1]{BeGe}. Before doing so we remark that the Lie algebra $\frg$ is assumed to be semisimple in \cite{BeGe}, whereas we consider the more general case of a reductive Lie algebra. The key result \cite[4.1]{BeGe}, however, remains true also for reductive Lie algebras, as we consider here $\frz$-finite modules. 

\vskip8pt

We will now explain how assumptions (i)-(iii) above are related to the assumptions (a), (b), (c) in \cite[4.1]{BeGe}. To connect the notation of \cite{BeGe} with that used here, we set $\theta = \chi_\lambda$ and $\theta' = \chi_\mu$. Note that the Verma module $M_{\lambda+\rho}$, as defined in \cite[1.9]{BeGe}, is equal to the Verma module $M(\lambda)$ with highest weight $\lambda$. Let $\eta$ be the map from $\frt^*_E$ to the set of characters of $\frz_E$, as introduced in \cite[1.6]{BeGe}.\footnote{Note that in \cite{BeGe} the untwisted Harish-Chandra homomorphism is used, in contrast to the twisted Harish-Chandra homomorphism used in \cite{Hu}, cf. \cite[1.9]{Hu}.} We then have $\eta(\lambda+\rho) = \theta$ and $\eta(\mu+\rho) = \theta'$. Put $\chi = w_\circ(\lambda+\rho)$ and $\psi = w_\circ(\mu+\rho)$, where $w_\circ$ is the longest Weyl group element. Since $\eta$ is invariant under the {\it linear} action of $W$, we also have $\eta(\chi) = \theta$ and $\eta(\psi) = \theta'$. Since $w_\circ(\Phi^+) = -\Phi^+$, and because the pairing $\pair{\cdot}{\cdot} $ is $W$-invariant, we find that $\pair{\chi}{\alpha^\vee} \notin \Z_{\le -1}$, and the same is true for $\psi$. In other words, $\chi$ and $\psi$ are dominant in the sense of \cite[1.5]{BeGe}.\footnote{Note that the concept of {\it dominant weights} in \cite{BeGe} differs by that defined in \cite[3.5]{Hu} by a shift by $\rho$.} Hence $\chi$ and $\psi$ satisfy condition (b) in \cite[4.1]{BeGe}. Furthermore, $\psi-\rho = w_\circ(\mu - \lambda)$ is integral because $\mu - \lambda$ is integral, and therefore condition (a) in \cite[4.1]{BeGe} is fulfilled. 

\vskip8pt

Note that $v \in W$ is contained in $W^\circ_\lambda$ if, by definition, $v \cdot \lambda = \lambda$, which is equivalent to $v(\lambda +\rho) = \lambda+\rho$, which is in turn equivalent to $w_\circ v w_\circ^{-1}$ being in the stabilizer $W_\chi$ of $\chi = w_\circ(\lambda+\rho)$ for the {\it linear} Weyl group action.\footnote{$W_\lambda$ being the notation introduced in \cite[1.4]{BeGe}.} Hence $W_\chi = w_\circ W^\circ_\lambda w_\circ^{-1}$. Condition (iii) above therefore implies (and is equivalent to) $W_\chi = W_\psi$, which is condition (c) in \cite[4.1]{BeGe}. 

\vskip8pt

\cite[4.1]{BeGe} then tells us that the functors above are equivalences of categories once we notice that $\mu - \lambda$ and $\psi - \chi$ are in the same linear $W$-orbit for the linear action, and hence determine the same finite-dimensional representation of $\frg_E$.  
\end{proof}

\vskip8pt

\subsection{Equivalences for categories of \texorpdfstring{$\DG$}{}-modules}
\begin{theorem}\label{D-equivalence}
Suppose $\lambda, \mu \in \frt^*_E$ satisfy conditions (i)-(iii) in \ref{BeGe theorem}. Then the functors 

\[\begin{array}{llcl}
\cTlamu: &\Dlambda & \lra &\Dmu \;,\\ 
&&&\\
\cTclamu: &\DGcoadlambda &\lra &\DGcoadlambda \;,\\ 
&&&\\
\cTrlamu: & \DrHfglambda & \lra &\DrHfgmu 
\end{array}\]

\vskip8pt

are equivalences of categories. One has natural isomorphisms $\cTlamu \circ \cTmula \cong \id$ and $\cTlamu \circ \cTlamu \cong \id$, and similarly for 
$\cTclamu$ and $\cTrlamu$.
\end{theorem}

\begin{proof}
The assertion for $\cTclamu$ follows immediately from the assertion for $\cTlamu$. We will first only treat the case of the functor $\cTlamu$.

\vskip8pt

By \ref{T adjoint}, $\cTlamu$ and $\cTmula$ are adjoints of each other, and so we have a natural transformation of functors $\vep: \cT := (\cTmula \circ \cTlamu) \ra \id$. In particular, for $M \in \Dlambda$ we obtain a morphism $\vep_M: \cT(M) \ra M$. 

\vskip8pt 

We consider the functors 

\[\begin{array}{cccc}
  \cB: \Uglambda & 
 \lra & \Dlambda\\
    N &\rightsquigarrow & \DG \otimes_{\UgE} N
\end{array}\]

\vskip8pt

For $M = \cB(N)$ we note that 

\[\begin{array}{ccccc}
\cT(M)= \cT(\cB(N)) & = & (\cTmula \circ \cTlamu)(\cB(N)) & \stackrel{\ref{T commute basechange}}{\cong} & (\cTmula \circ \cB \circ \Tlamu) (N)   \\
&&&&\\
& \stackrel{\ref{T commute basechange}}{\cong} & (\cB \circ  \Tmula \circ \Tlamu) (N) & \stackrel{\ref{BeGe theorem}}{\cong} & \cB(N) = M \;.
\end{array}\]

\vskip8pt

Given an arbitrary $M \in \Dlambda$, we denote by $N$ the underlying $\UgE$-module. Consider the canonical map 

\[\phi: M_1 := \cB(N) = \DG \ot_{\UgE} N \lra M\,, \;\, \delta \ot n \mapsto \delta.n\]

\vskip8pt

of $\DG$-modules, which is evidently surjective. Let $M_2$ be its kernel, so  we have an exact sequence

\[0 \ra M_2 \ra M_1 \ra M \ra 0\]

\vskip8pt

of $\DG$-modules. By \ref{translation exact} we obtain a commutative diagram with exact rows 

\[\begin{tikzcd}
0  \arrow{r} & \cT(M_2) \arrow{r}\arrow{d}{\vep_{M_2}}  &  \cT(M_1)  \arrow{r}{\cT(\phi)}\arrow{d}{\vep_{M_1}} &    \cT(M) \arrow{r}\arrow{d}{\vep_{M}} & 0 \\
 0  \arrow{r} & M_2 \arrow{r}  &  M_1 \arrow{r}{\phi} &    M \arrow{r} & 0 \;\;.
\end{tikzcd}\]

\vskip8pt

The discussion above shows that $\vep_{M_1}$ is an isomorphism. 

\vskip8pt

Let $m \in M$. Since $\phi$ is surjective there is a $m_1 \in M_1$ such that $\phi(m_1)= m$. Let $m_1' \in \cT(M_1)$ be such that $\vep_{M_1}(m_1')= m_1$. Then by commutativity of the diagram we see that $\vep_M(\cT(\phi)(m_1')) = m$. Thus $\vep_M$ is surjective. Therefore, $\vep_M$ is surjective for {\it every} $M \in \Dlambda$. In particular,  $\vep_{M_2}$ is surjective. By four lemma we get that $\vep_M$ is injective as well. Thus $\vep_M$ is an isomorphism for any $M \in \Dlambda$.

\vskip8pt

Exactly the same arguments as above show that 

\[\cTrlamu: \DrHlambda \lra \DrHmu\]

\vskip8pt

is an equivalence of categories with quasi-inverse $\cTrmula$. Because $\cTrlamu$ is an endo-functor of $\DrHfgzfin$, by \ref{zfinite stable}, it induces an equivalence $\DrHfglambda \ra \DrHfgmu$ with quasi-inverse $\cTrmula$.  
\end{proof}

\vskip8pt

\section{Relation with translation functors on category \texorpdfstring{$\cO$}{}}

\subsection{Category \texorpdfstring{$\cOPinfty$}{} and the functor \texorpdfstring{$\chcFGP$}{}}
Let $\bP \sub \bG$ be the standard parabolic subgroup with Lie algebra $\frp$. 

\begin{para}\label{cOpinfty}{\it Category $\cOpinfty$ for the reductive Lie algebra $\frg$.}
Given a representation $\phi: \frt_E \ra \End_E(M)$, a weight $\lambda \in \frt^*_E$, and a positive integer $i \ge 1$ we set 

\[M^i_\lambda = \{m \in M \midc \forall \frx \in \frt_E: (\phi(\frx)-\lambda(\frx)\cdot \id)^i.m = 0\} \;.\]

\vskip8pt

$M^\infty_\lambda := \bigcup_{i \ge 1} M^i_\lambda$ is then the generalized eigenspace for the weight $\lambda$.

\vskip8pt

The category $\cOpinfty$ for the pair $(\frg,\frp)$ and the coefficient field $E$ is defined to be the full subcategory of all $\UgE$-modules $M$ which satisfy the following properties:

\vskip8pt

\begin{enumerate}
    \item  $M$ is finitely generated as a $\UgE$-module.
    \item  $M = \bigoplus_{\lambda \in \frt^*_E} M_\lambda^\infty$ \;. 
    \item  The action of $\frp_E$ on $M$ is locally finite, i.e. for every $m\in M$, the subspace $\UpE .m \subset M$ is finite-dimensional over $E$.
\end{enumerate}

\vskip8pt

$\cOpinfty$ contains the full subcategory $\cOp$ consisting of all modules $M$ satisfying the condition

\[\hskip-20pt (2)^1 \hskip20pt M = \bigoplus_{\lambda \in \frt^*_E} M_\lambda^1\] 

\vskip8pt

instead of (2) above. It is shown in \cite[3.2.4]{AS} that an $\UgE$-module $M$ is contained in $\cOpinfty$ if and only if $M$ has a finite filtration whose successive subquotients are in $\cOp$. This implies that $\cOpinfty$ is stable under tensoring with finite-dimensional $\UgE$-modules and that it is a subcategory of $\Ugzfin$ (cf. \cite[1.1]{Hu} for the case of category $\cO$). Category $\cOpinfty$ is a also Serre subcategory of $\UgEmod$ \cite[3.1.2]{AS}. 

\vskip8pt

Furthermore, $\cOpinfty$ contains the full subcategory $\cOpinftyalg$ consisting of all modules $M$ satisfying the stronger condition

\[\hskip-20pt (2)_\alg \hskip20pt M = \bigoplus_{\lambda \in \frt^*_E \; {\rm algebraic}} M_\lambda^\infty\] 

\vskip8pt

instead of (2) above. $\cOpinftyalg$ is also a Serre subcategory of $\UgEmod$ \cite[3.1.4]{AS}. 
\end{para}

\vskip8pt

\begin{para}{\it Category $\cOPinfty$.}
Here we follow ideas of \cite{O-S3} but place ourselves in the more general setting of category $\cOpinfty$. We denote by $\cOPinfty$ the category of pairs $(M,\phi)$, where $M$ is an object of $\cOpinfty$ and $\phi: P = \bP(F) \ra \End_E(M)^\x$ is a group homomorphism satisfying the following conditions:

\vskip8pt

\begin{enumerate}
    \item[(i)] $M$ possesses an exhaustive ascending filtration $w_\circ 
    \sub W_1 \sub W_2 \sub \ldots$ by finite-dimensional subspaces $W_j$ which are stable under the action of $P$ via $\phi$;
    \item[(ii)] for each $j \ge 1$ the action of $P$ on $W_j$ is locally $F$-analytic;
    \item[(iii)] the derived action ${\rm d} \phi: \frp \ra \End_E(M)$ coincides with the restriction of the $\frg_E$-module structure of $M$ to $\frp$;
    \item[(iv)] $\forall g \in P, \, \frx \in \frg_E, \, m \in M: \; \phi(h).(\frx.m) = \Ad(h)(\frx).(\phi(h).m)$.
\end{enumerate}

\vskip8pt

A morphism $f: (M,\phi_M) \ra (N,\phi_N)$ in $\cOPinfty$ is a morphism $M \ra N$ of the underlying $\UgE$-modules which commutes with the actions of $P$ on $M$ and $N$. By \cite[7.4.2]{AS}, any object in $\cOPinfty$ carries the structure of a $\DgP$-module, which extends the module structures over $\UgE$ and $E[P]$. Considering an object of $\cOPinfty$ as a $\DgP$-module gives a fully faithful functor $\cOPinfty \ra \DgP\mbox{-}{\rm mod}$.
\end{para}

\vskip8pt

\begin{example}
Let $\log: T \ra \frt_E$ be a logarithm on $T$ as defined in \cite[2.2.1]{AS}. Generalizing the work done in \cite{O-S2} for category $\cO^\frp_\alg$, it was shown in \cite[2.4.16]{AS} that, depending on the choice of $\log$, every $M$ in $\cOpinftyalg$ has a natural lift $\Lift(M,\log)$ which is an object of $\cOPinfty$, and $\Lift(-,\log): \cOpinftyalg \ra \cOPinfty$ is a functor.   
\end{example}

\vskip8pt

\begin{example}\label{tau Vermas}
Let $\tau: P \ra \Ex$ be a locally $F$-analytic character. Let $\rmd \tau: \frp \ra E$ be the derivative of $\tau$, and denote by $E_{\rmd \tau}$ the corresponding 1-dimensional $\UpE$-module. Set $M = \UgE \otimes_{\UpE} E_{\rmd \tau}$ (a generalized Verma module), and define a $P$-action on $M$ by $g.(u \ot 1) = \Ad(g)(u) \ot \tau(g)$, for $g \in P$ and $u \in \UgE$. Together with this $P$-action $M$ becomes an object of $\cOPinfty$ which we denote by $M^P(\tau)$. This object does not lie in the image of the previous functor $\cOpinftyalg \ra \cOPinfty$, unless $\tau|_T$ is algebraic.  
\end{example}

\vskip8pt

\begin{para}\label{tensoring for cOPinfty}{\it Translation functors for category $\cOPinfty$.} Let $L$ be a finite-dimensional locally $F$-analytic representation of $G$. Given a  $M$ in $\cOPinfty$ we equip $L \ot_E M$ with the diagonal action of $P$. Because $\cOpinfty$ is stable under tensoring with finite-dimensional modules, it follows that $\cOPinfty$ is stable under tensoring with finite-dimensional locally $F$-analytic representations of $G$. In particular, for compatible $\lambda, \mu \in \frt^*_E$ the functor $\Tlamu$ from \ref{Tlamu} maps $\cO^{P,\infty}_{|\lambda|}$ to $\cO^{P,\infty}_{|\mu|}$.  
\end{para}

\begin{para}{\it The functor $\chcFGP$.} Recall that a locally $F$-analytic representation $V$ of $G$ is called {\it strongly admissible}, if $V$ is of compact type and $V'_b$ is finitely generated as module over $\DH$ for any (equivalently, one) compact open subgroup $H \sub G$, cf. \cite[p. 453-454]{S-T1}. Denote by $\Rep^{\sm}_E(L_P)^\sadm$ the category of smooth strongly admissible representations on $E$-vector spaces of the Levi subgroup $L_P$ of $P$. In the spirit of \cite{O-S2,O-S3,AS} we introduce the functor 

\[\begin{array}{ccc}
 \chcFGP   : \cOPinfty \times \Rep^{\sm}_E(L_P)^\sadm & \lra & \DGmod  \\
     (M,V) & \rightsquigarrow & \DG \ot_{\DgP} (M \ot_E V') \;.
\end{array}\]

\vskip8pt

When $V$ is the trivial one-dimensional representation $\triv$ of $L_P$, we write $\chcFGP(M)$ instead of $\chcFGP(M,\triv)$. 
\end{para}

\vskip8pt

\begin{example} Let $\tau: P \ra \Ex$ be a locally $F$-analytic character. $\tau$ gives rise to a morphism of $E$-algebras $\DP \ra E$, which sends $g \in P$ to $\tau(g)$. We denote by $E_\tau$. It follows from \ref{U-DP-bimodule iso} that the inclusion $\UgE \hra \DgP$ induces an isomorphism

\begin{numequation}\label{Verma iso}
\UgE \ot_{\UpE} E_{\rmd \tau} \lra \DgP \ot_{\DP} E_\tau
\end{numequation}

\vskip8pt

of $\UgE$-modules. In fact, it is straightforward to see that \ref{Verma iso} is in fact also an isomorphism of $P$-representations, when we give the left hand side the $P$-module structure as in \ref{tau Vermas}. Hence $M^P(\tau) = \DgP \ot_{\DP} E_\tau$. This implies

\begin{numequation}\label{FGP on Verma}
\begin{array}{rcl} \chcFGP(M^P(\tau)) & = & \DG \ot_{\DgP} M^P(\tau)\\
&&\\
& \stackrel{\ref{Verma iso}}{=}& \DG \ot_{\DgP} (\DgP \ot_{\DP} E_\tau)\\
&&\\
& = & \DG \ot_{\DP} E_\tau \;.
\end{array}
\end{numequation}
\end{example}

\vskip8pt

\begin{prop}
For  any $(M,V)$ in $\cOPinfty \times \Rep^{\sm}_E(L_P)^\sadm$ the module $\chcFGP(M,V)$ belongs to $\DGcoadzfin$.
\end{prop}

\begin{proof} As we have noticed in \ref{cOpinfty}, objects in $\cOpinfty$ are $\frz$-finite. Note that $\UgE$ acts via the augmentation $\UgE \ra E$ on $V'$, and the same is true for $\frz_E$. Hence $M \ot_E V'$ is $\frz$-finite. Since $\frz_E \sub Z(\DG)$, hence $\frz_E \sub Z(\DgP)$, it follows that  $\chcFGP(M,V)$ is $\frz$-finite.

\vskip8pt

If $M = \Lift(N,\log)$ for $N$ in $\cOpinftyalg$, then the assertion that $\chcFGP(M,V)$ is coadmissible, is \cite[4.2.3]{AS}. We note that the proof of the key result \cite[4.1.5]{AS} actually does not use that $M$ is of this particular form and applies more generally to all $M$ in $\cOPinfty$.
\end{proof}

\vskip8pt

\begin{rem}
While we will mostly work with the functor $\chcFGP$, the previous result, together with the categorical anti-equivalence of \cite[6.3]{S-T2} allows us to pass to continuous dual spaces (equipped with the strong topology) and thus consider the functor 

\[\begin{array}{ccc}
   \cF^G_P: \cOPinfty \times \Rep^{\sm}_E(L_P)^\sadm & \lra &\Rep^\la_E(G)^{\rm coad}\\
     (M,V) & \rightsquigarrow & \Big(\DG \ot_{\DgP} (M \ot_E V')\Big)'_b
\end{array}\]

\vskip8pt

which takes on values in the category of admissible locally $F$-analytic representations. 
\qed\end{rem}

\vskip8pt

\begin{theorem}\label{F-T commute}
Let $\lambda, \mu \in \frt^*_E$ be compatible. For any $(M,V)$ in $\cOPinfty \times \Rep^{\sm}_E(L_P)^\sadm$ there is a canonical isomorphism

\[ \cTlamu\Big(\chcFGP(M,V)\Big) \cong \chcFGP \Big(T_\lambda^\mu(M),V\Big)\]

\vskip8pt

which is natural in $M$ and $V$.
\end{theorem}

\vskip8pt

\begin{proof}
Set $N = M \ot_E V'$. By definition of $\Tlamu$ and $\chcFGP$ we have

\[\begin{array}{rcl}
\chcFGP \Big(\Tlamu(M),V\Big) & = & \DG\ot_{\DgP} \pr_{|\mu|}\Big(L(\onu)\ot_E \pr_{|\lambda|}(M)\Big) \ot_E V'\\
&&\\
& = &  \DG\ot_{\DgP} \pr_{|\mu|}\Big(L(\onu)\ot_E \pr_{|\lambda|}(M) \ot_E V'\Big)\\
&&\\
& = &  \DG\ot_{\DgP} \pr_{|\mu|}\Big(L(\onu)\ot_E \pr_{|\lambda|}(M\ot_E V')\Big)\\
&&\\
& \stackrel{\ref{Homzeroes}}{\cong} & \pr_{|\mu|} \Big(\DG\ot_{\DgP}(L(\onu)\ot_E \pr_{|\lambda|}(N)\Big) \\
&&\\
& \stackrel{\ref{finite commuting}}{\cong} & \pr_{|\mu|}\Big(L(\onu) \ot_E (\DG\ot_{\DgP}\pr_{|\lambda|} (N))\Big)\\
&&\\
& \stackrel{\ref{Homzeroes}}{\cong} & \pr_{|\mu|}\Big(L(\onu) \ot_E \pr_{|\lambda|} (\DG \ot_{\DgP} N)\Big)\\
&&\\
& = & \cTlamu\Big(\chcFGP(M,V)\Big)
\end{array}\]

\vskip8pt

where the last equality holds by the definition of $\cTlamu$.
\end{proof}

\vskip8pt

\subsection{Locally analytic principal series representations}

In this section we study more closely the effect of translation functors on  locally analytic principal series representations, or rather their associated coadmissible modules.

\begin{para}\label{inducedrep} For a locally $F$-analytic character $\tau: T \ra E^\times$ we consider the locally analytic principal series representation

\[\Ind^G_B(\tau) = \Big\{f: G \ra E \mbox{ locally $F$-analytic} \, \Big |\, \forall \, g \in G \; \forall\; b \in B: \; f(bg) = \tau(b)f(g)\Big\} \,.\]

\vskip8pt

The action of $G$ on $\Ind^G_B(\tau)$ is given by $g.f(x)= f(x. g^{-1})$. It is an easy consequence of $B\bksl G$ being compact that $\Ind^G_B(\tau)$ is strongly admissible and hence admissible. For the corresponding coadmissible module there is a canonical isomorphism of $\DG$-modules

\begin{numequation}\label{Tensor description of induction}
\Ind^G_B(\tau)'_b \cong \DG \ot_{\DP} E_{\tau^{-1}}\;,
\end{numequation}

cf. \cite[6.1 (iv)]{ST_duality}.
\end{para}

\vskip8pt

\begin{para}\label{dot action loc an}{\it The dot-action of the Weyl group on characters.}
Following standard conventions we set $\rho = \frac{1}{2}\sum_{\alpha \in \Phi^+} \alpha$. It is well-known and easy to prove that for any $w \in W$ one has $w(\rho)-\rho \in \Lambda_r$. Since weights in $\Lambda_r$ are algebraic, there is a unique algebraic character $\gamma_w \in X^*(\bT)$ such that $\rmd \gamma_w = w(\rho)-\rho$. Denote by $v_F$ the $p$-adic valuation on $F$ normalized by $v_p(\Fx) = \Z$. Choose a uniformizer $\vpi_F$ of $F$ and set $|x| = \vpi_F^{-v_F(x)}$ for all $x \in \Fx$, so that $v_F(x|x|) = 0$ for all $x \in \Fx$. When $F = \Qp$ we choose $\vpi_\Qp = p$. Denote by $|\gamma_w|$ the smooth character $T \ra \Ex$, $t \mapsto |\gamma_w(t)|$. 

\vskip8pt

Let $\tau: T \ra \Ex$ be a locally $F$-analytic character of $T$. We define the dot action of the Weyl group $W$ on $\tau$ as 

\begin{numequation}\label{dot action}
w\cdot \tau = w(\tau) \gamma_w |\gamma_w| \;,
\end{numequation}

\vskip8pt

where $w(\tau)(t) = \tau(\dot{w}^{-1}t\dot{w})$, for any $\dot{w} \in N_G(T)$ representing $w$. Note that ${\rm d}(w\cdot \tau) = w\cdot({\rm d}\tau)$, so the dot-action on locally $F$-analytic characters and on $\frt^*_E$ are compatible via the derivative.
\end{para}

\vskip8pt 

\begin{rem} As $|\gamma_w|$ is a smooth character, the compatibility with the dot-action on weights would also be given if we had defined $w \cdot \tau$ as $w(\tau)\gamma_w$. To define the dot-action on locally analytic characters by the formula \ref{dot action} (with the factor $|\gamma_w|$) has been suggested by the relations to 2-dimensional Galois representations in sec. \ref{Galois}. Moreover, with the formula \ref{dot action}, unitary characters are preserved under the dot-action.
\qed\end{rem}

\vskip8pt

\begin{para}{\it The Weyl groups $\Wbrla$, their associated Euclidean spaces, and facets.} Recall the Euclidean vector space $\sE = \bbR \otimes_\Z \Phi$. Following \cite[3.4]{Hu} we set for $\lambda \in \frt^*_E$

\[\Phibrla = \{\alpha \in \Phi \midc \pair{\lambda}{\alpha^\vee} \in \Z\} \;, \; \Wbrla = \{w\in W \midc w \cdot \lambda - \lambda \in \Lambda_r\} \;, \; \sE(\lambda) = \langle \Phibrla \rangle_\bbR \;.\]

\vskip8pt

If $\mu-\lambda \in \Lambda$, in particular if $\lambda$ and $\mu$ are compatible, then $\Wbrla = W_{[\mu]}$. It is known that $\Phibrla$ is a root system in its $\bbR$-span $\sE(\lambda) \sub \sE$ with Weyl group $\Wbrla$, and $\Phibrla^+ = \Phi^+ \cap \Phibrla$ is a system of positive roots in $\Phibrla$ \cite[3.4, Thm. and Exercise]{Hu}. A non-empty subset $\sF$ of $\sE(\lambda)$ is called a {\it facet} if there is a partition $\Phibrla^+ = \Phi^0_{[\lambda],\sF} \sqcup \Phi^+_{[\lambda],\sF} \sqcup \Phi^-_{[\lambda],\sF}$ such that 

\[\lambda \in \sF \; \Longleftrightarrow \; \left\{\begin{array}{lcl} \pair{\lambda+\rho}{\alpha^\vee} = 0 & \mbox{ when } & \alpha \in \Phi^0_{[\lambda],\sF} \\
\pair{\lambda+\rho}{\alpha^\vee} > 0 & \mbox{ when } & \alpha \in \Phi^+_{[\lambda],\sF} \\
\pair{\lambda+\rho}{\alpha^\vee} < 0 & \mbox{ when } & \alpha \in \Phi^-_{[\lambda],\sF} \\
\end{array} \right.\]

\vskip8pt

$\sE(\lambda)$ is the disjoint union of its facets, and $\Wbrla$ maps facts to facets under the dot-action. There is a unique $\lambda^\natural \in \sE(\lambda)$ such that 

\[\pair{\lambda^\natural}{\alpha^\vee} = \pair{\lambda}{\alpha^\vee} \;.\]
\end{para}

\vskip8pt

The following key lemma from \cite[7.5]{Hu} will also be important for our purposes here.

\begin{lemma}\label{7.5}
Let $\lambda, \mu\in \frt^*_E$ be compatible in the sense of \ref{compatible}. Set $\nu  = \mu - \lambda$, and let $\onu$ be the unique conjugate of $\nu$ in $\Lambda^+$ under the linear $W$-action. Assume that 

\begin{numequation}\label{key cond}
\mbox{$\mu^\natural$ lies in the closure of the facet $\sF \sub \sE(\lambda)$ which contains  $\lambda^\natural$ \;.}
\end{numequation}

Then for all weights $\nu_1 \neq \nu$ of $L(\onu)$ the weight $\nu_1 + \lambda$ is not conjugate to $\nu + \lambda = \mu$ by the dot-action of $W_{[\lambda]} = W_{[\mu]}$.
\end{lemma}

\vskip8pt

\begin{cor}\label{cor7.5} With the notation and under the assumptions of \ref{7.5}, for any $w \in \Wbrla$ and any weight $\nu_1 \neq w(\nu)$ of $L(\onu)$, the weight $\nu_1 + w \cdot\lambda$ is not conjugate to $\mu$ under the dot-action of the whole Weyl group $W$. 
\end{cor}

\begin{proof}
Because $(w \cdot \lambda)^\natural = w \cdot \lambda^\natural$ and $(w \cdot \mu)^\natural = w \cdot \mu^\natural$ \cite[7.4]{Hu}, we find that $(w \cdot \mu)^\natural$ is contained in the closure of the facet $w \cdot \sF$ to which $(w \cdot \lambda)^\natural$ belongs. We can then apply \ref{7.5} to $w \cdot \lambda$, $w \cdot \mu$ and $w(\nu) = w \cdot \mu - w \cdot \lambda$, and deduce that for $\nu_1 \neq w(\nu)$ the weight $\nu_1 + w\cdot \lambda$ is not conjugate to $w\cdot \mu$ under the dot-action of $\Wbrla$. But then it is not conjugate to $\mu$ under the dot-action of $\Wbrla$.

\vskip8pt

Now assume $u \in W$ is such that $u \cdot (\nu_1 + w \cdot \lambda) = \mu$. This is equivalent to $(uw) \cdot \lambda - \lambda = \nu - u(\nu_1)$. But since $\nu$ and $u(\nu_1)$ are both weights of $L(\onu)$, and because $L(\onu)$ is a highest weight module, $\nu - u(\nu_1) \in \Lambda_r$. Therefore $uw \in \Wbrla$ and hence $u \in \Wbrla$. Hence there is no such $u \in W$ unless $\nu_1 = w(\nu)$.
\end{proof}

\vskip8pt

\begin{prop}\label{T on tau Vermas} Let $\wlambda, \wmu: T \ra \Ex$ be locally $F$-analytic characters of $T$, and set $\lambda = \rmd \wlambda$, $\mu  = \rmd \wmu$. Suppose that $\wnu := \wmu \wlambda^{-1}$ is algebraic and that $\lambda$ and $\mu$ are both anti-dominant. Assume that condition \ref{key cond} is fulfilled. Then, for all $w \in W_{[\lambda]} = W_{[\mu]}$ there is an isomorphism of $\DgB$-modules

\[\Tlamu\Big(M^B(w \cdot \wlambda)\Big) \cong M^B(w \cdot \wmu) \;.\]
\end{prop}

\vskip8pt

\begin{proof} We follow the arguments given in the proof of \cite[7.6]{Hu}. The finite-dimensional representation $L(\onu)$ has a $B$-stable filtration 

\begin{numequation}\label{filtration L}
0 = L_0 \subsetneq L_1 \subsetneq \ldots \subsetneq L_d = L(\nu)
\end{numequation}

such that for any $i \in \{1, \ldots,d\}$ the quotient $L_i/L_{i-1}$ is one-dimensional and the action of $B$ factors through $T$ and is given by $\wnu_1$, where $\nu_1$ runs through the weights of $L(\onu)$, $\wnu_1$ is the (unique) lift of $\nu_1$ to an algebraic character of $\bT$, and this quotient appears with multiplicity $\dim_E(L(\onu)_{\nu_1})$. We use the canonical isomorphism of $\DgB$-modules from \ref{finite commuting}

\[L(\onu) \otimes_E M^B(w \cdot \wlambda) \cong \DgB \ot_{\DB} \Big(L(\onu)|_B \ot_E E_{w \cdot \wlambda}\Big) \;.\]

\vskip8pt

Using the filtration \ref{filtration L} we see that $L(\onu) \otimes_E M^B(w \cdot \wlambda)$ has a filtration of $\DgB$-submodules $M_i = \DgB \ot_{\DB} (L_i \ot_E E_{w \cdot \wlambda})$ such that 

\[M_i/M_{i-1} \cong M^B\Big(\wnu_1 (w \cdot \wlambda)\Big)\]

\vskip8pt

for some weight $\nu_1$ of $L(\onu)$. As a $\UgE$-module $M^B\Big(\wnu_1 (w \cdot \wlambda)\Big)$ is equal to $M(\nu_1+ w\cdot \lambda))$. For $M(\nu_1+ w\cdot \lambda))$ to lie in the subcategory 
$\Ugmod_{|\mu|}$ it is necessary and sufficient that $\nu_1+ w\cdot \lambda$ is conjugate to $\mu$ under the dot-action of $W$. By \ref{cor7.5} this is only the case when $\nu_1 = w(\nu)$. Therefore,

\[\pr_{|\mu|}\Big(L(\onu) \otimes_E M^B(w \cdot \wlambda)\Big) \cong M^B\Big(w(\wnu) (w \cdot \wlambda)\Big) = M^B(w \cdot \wmu) \;.\] 
\end{proof}

\vskip8pt

\begin{theorem}\label{translation of induced representations}
Let $\wlambda, \wmu: T \ra \Ex$ be locally $F$-analytic characters of $T$, and set $\lambda = \rmd \wlambda$, $\mu  = \rmd \wmu$. Suppose $\wnu := \wmu \wlambda^{-1}$ is algebraic and $\lambda$ and $\mu$ are both anti-dominant. Assume that condition \ref{key cond} is fulfilled. Then, for all $w \in W_{[\lambda]} = W_{[\mu]}$ there is an isomorphism of $\DG$-modules

\[\cTlamu\Big(\DG\ot_{D(B)} E_{w \cdot \wlambda}\Big) \cong \DG\ot_{D(B)} E_{w \cdot \wmu} \;.\]
\end{theorem}

\vskip8pt

\begin{proof} 
This follows from \ref{FGP on Verma}, \ref{F-T commute}, and \ref{T on tau Vermas}.
\end{proof}

\section{\texorpdfstring{The case $\GL_2(\Qp)$: relations with Galois representations}{}}\label{Galois}

We denote by $\sG_\Qp$ the absolute Galois group of $\Qp$. For a 2-dimensional absolutely irreducible representation $V$ over $E$ of $\sG_\Qp$ let $\Pi(V)$ be the corresponding unitary Banach space representation of $\GL_2(\Qp)$ associated to $V$ via the $p$-adic local Langlands correspondence of \cite{Colmez_GL2,ColDoPa}. We let $\Pi(V)^\la$ be the associated locally analytic representation which is strongly admissible. Our aim in this section is to study the effect of the translation functors on $\Pi(V)^\la$ or rather on the corresponding coadmissible module.

\vskip8pt

\subsection{2-dimensional trianguline representations of \texorpdfstring{$\sG_\Qp$}{}}

\begin{para}\label{Trianguline variety of Colmez}{\it The trianguline variety.}
Following Colmez \cite{Co1} let $\hsT(E)$ be the set of continuous characters\footnote{Any continuous character on $\Qpx$ is locally analytic by \cite[9.4]{DDMS}.} $\delta: \Qpx \ra E^{*}$. $\hsT(E)$ is the set of $E$-rational points of a rigid analytic variety $\hsT$ over $\Qp$. We denote by $x$ the character given by $\Qpx \hra \Ex$ and by $|x|$ the character $x\ra p^{-v_p(x)}$.

\vskip8pt

In the following we denote by $(\delta_1, \delta_2)$ elements of $\hsT(E) \times \hsT(E)$. Define a subset $\sZ := \sZ_\sp \sqcup \sZ'_\sp \sub \hsT(E) \times \hsT(E)$ by

\[\begin{array}{rcl}
\sZ_\sp & := & \bigcup_{i \in \Z_{\ge 1}} \{(\delta_1,\delta_2) \midc \delta_1. \delta_2^{-1} = x^i |x| \} \;,\\
&&\\
\sZ'_\sp & := &\bigcup_{i \in \Z_{\ge 0}} \{(\delta_1,\delta_2) \midc \delta_1. \delta_2^{-1} = x^{-i}\} \;.
\end{array}\]

\vskip8pt

Set 

\[\sS := \Big(\hsT(E) \times \hsT(E) \setminus \sZ\Big) \times \{\infty\} \; \sqcup \; \sZ \times \bbP^1(E) \;.\]

\vskip8pt

Let $\pi: \sS \ra \hsT(E) \times \hsT(E)$ be the map given by $\pi(\delta_1,\delta_2,\sL) = (\delta_1,\delta_2)$.

\vskip8pt
 
For a locally analytic character $\delta: \Qpx \ra \Ex$ we set $w(\delta) := \frac{\log(\delta(u))}{\log (u)}$, where $u\in \bbZ_p^*$ and is not a root of unity.\footnote{The use of $w$ here should not lead to a confusion with the use of $w$ for elements of Weyl groups. Here there is only one non-trivial Weyl group element, and we denote it by $w_\circ$.} This is the same as the derivative of $\delta$ at 1. Set 

\[\sS_{*} :=\{(\delta_1,\delta_2,\sL) \in \sS \midc v_p(\delta_1(p))+ v_p(\delta_2(p)) = 0 \,,\, v_p(\delta_1(p)) > 0\} \;.\]
 
\vskip8pt
 
For $s = (\delta_1, \delta_2, \sL)\in \sS_{*}$ put $u(s):= v_p(\delta_1(p))$ and $w(s):= w(\delta_1)- w(\delta_2)$. Then define
 
\[\begin{array}{lcl}
\sS_{*}^{\rm ng} & := & \{s\in \sS_{*} \midc \,w(s) \notin \bZ_{\ge 1}\}  \\
&&\\      
\sS_{*}^{\rm cris} & := & \{s\in \sS_{*}\,|\, w(s) \in \bZ_{\ge 1}, u(s) < w(s), \; \sL=\infty\}\\
&&\\      
\sS_{*}^{\rm st} & := & \{s\in \sS_{*} \midc w(s) \in \bZ_{\ge 1}, u(s) < w(s), \; \sL\neq\infty\}\\
&&\\      
\sS_{\irr} & := & \sS_{*}^{\rm ng} \; \bigsqcup \; \sS_{*}^{\rm cris} \; \bigsqcup \; \sS_{*}^{\rm st}
\end{array}\]
\end{para}

\vskip8pt

\begin{rem}\label{no *-pts over Z'} One has $\sS_* \cap \pi^{-1}(\sZ'_\sp) = \emptyset$. Indeed, suppose $s = (\delta_1,\delta_2,\sL) \in \sS_*$ and 
$\delta_2 = x^i\delta_1$ for some integer $i \ge 0$. Then $v_p(\delta(p)) = v_p(\delta_1(p) + i$, hence $0 = 2v_p(\delta_1(p)) + i$ and thus $v_p(\delta(p)) = -\frac{i}{2} \le 0$, contradicting the second condition for $s$ to belong to $\sS_*$.  
\qed\end{rem}

\vskip8pt

\begin{para}\label{generic}
We let $\bG = \GL_{2,\Qp}$, $\bB$ the Borel subgroup of upper triangular matrices, and $\bT \sub \bB$ the torus of diagonal matrices. We write locally analytic characters $\vtheta$ of $T$ as $\vtheta = \vtheta_1 \otimes \vtheta_2$ if $\vtheta(\diag(a,d)) = \vtheta_1(a)\vtheta_2(d)$. For $(\delta_1,\delta_2) \in \hsT(E) \times \hsT(E)$ set 

\[B^{\la}(\delta_1, \delta_2):= \Ind^{G}_B(\delta_2\ot \delta_1(x|x|)^{-1}) \;.\]

\vskip8pt

In the following we set 

\begin{numequation}\label{lambdatilde}
\wlambda = (\delta_2\ot \delta_1(x|x|)^{-1})^{-1} = \delta_2^{-1} \ot \delta_1^{-1} x|x|
\end{numequation}

\vskip8pt

so that $B^{\la}(\delta_1, \delta_2) = \Ind^G_B(\wlambda^{-1})$. Let $w_\circ$ be the non-trivial Weyl group element. Recall that in \ref{dot action} we defined a dot-action of the Weyl group on the set of locally analytic characters. In general, it depends on the choice of the absolute value $|\cdot|$ on $F$ with values in $E$, normalized such that $|p| = p^{-1}$, but since here we have $F = \Qp$, the normalization determines the absolute value. Given a locally analytic character $\vtheta = \vtheta_1 \otimes \vtheta_2$ of $T$, we have 

\begin{numequation}\label{w0 action}
w_\circ \cdot \vtheta = \vtheta_2 (x|x|)^{-1} \otimes \vtheta_1 x|x| \;. 
\end{numequation}

To $s \in \sS_{\irr}$ there is attached an absolutely irreducible 2-dimensional representation $V(s)$ of $\sG_\Qp$ over $E$ \cite[0.5 (ii)]{Co1}.\footnote{It is not said there that $V(s)$ is absolutely irreducible, but this follows because the assertion is independent of the coefficient field.} We write $\Pi(s) := \Pi(V(s))$ for the absolutely irreducible unitary Banach space representation of $G$.   
\end{para}

\vskip8pt

$s = (\delta_1,\delta_2,\sL) \in \sS_{\irr}$ is called {\it generic} if $\pi(s) \notin \sZ_\sp$. By \ref{no *-pts over Z'}, generic $s$ are contained in $\sS_*^{\rm ng} \sqcup \sS_*^{\rm cris}$, and hence $\sL = \infty$ for those $s$. We have the following

\vskip8pt

\begin{theorem}\label{colmez 8.7} Let $s = (\delta_1,\delta_2,\infty) \in \sS_{\irr}$ be generic and $\wlambda$ as in \ref{lambdatilde}. 

\vskip8pt

(i) There is a non-split exact sequence of admissible locally analytic representations

\[0 \lra B^{\la}(\delta_1, \delta_2) \lra \Pi(s)^{\la} \lra B^{\la}(\delta_2, \delta_1) \lra 0 \;.\]

\vskip8pt

(ii) There is a non-split exact sequence of coadmissible $\DG$-modules

\begin{numequation}\label{colmez 8.7 seq}
0 \ra \DG \ot_{\DB} E_{w_\circ \cdot \wlambda} \ra (\Pi(s)^{\la})'_b \ra \DG \ot_{\DB} E_{\wlambda} \ra 0\;.
\end{numequation}

The Ext group $\Ext^1_{\DG}\Big(\DG \ot_{\DB} E_{\wlambda},\DG \ot_{\DB} E_{w_\circ \cdot \wlambda}\Big)$, computed in the category of abstract $\DG$-modules, is one-dimensional over $E$.
\end{theorem}

\vskip8pt

\begin{proof} (i) This is \cite[Thm 8.7 (i)]{Co2}.

\vskip8pt

(ii) The first assertion follows from (i) by passing to strong dual spaces, together with \ref{Tensor description of induction} and \ref{w0 action}. The second assertion follows from \cite[8.18, 8.15]{KohlhaaseCohomology}, cf. also \cite[8.8. (i)]{Co2}.
\end{proof}

\vskip8pt

\subsection{The effect of translation functors on \texorpdfstring{$\Pi(s)^\la$}{}}

\begin{para}\label{setup}{\it The set-up.} Let $w_1,w_2 \in \Z$ be integers and consider the algebraic character $\wnu = x^{w_1} \otimes x^{w_2}$ of $\bT$, so that $\nu := \rmd \wnu = w_1 \vep_1 + w_2 \vep_2$, where $\vep_1, \vep_2 \in \frt^*$ are defined by $\vep_1((\diag(a,d)) = a$, $\vep_2(\diag(a,d)) = d$. 

\vskip8pt

Let $s = (\delta_1,\delta_2,\infty) \in \sS_{\irr}$ be a generic point and $\wlambda = \delta_2^{-1} \ot \delta_1^{-1}x|x|$ be as in \ref{lambdatilde}. Set $\lambda  = \rmd \wlambda$ and $\mu = \nu + \lambda$. Our aim is to understand the action of $\cTlamu$ on the coadmissible module of $\Pi(s)^\la$. However, unless $w_1+w_2 = 0$, the center of $G$ does not act by a unitary character on $\cTlamu((\Pi(s)^\la)')$. But this is easy to achieve by twisting with an unramified smooth character of $G$, after possibly passing to a quadratic extension of $E$.\footnote{This is not the only way to obtain a module on which the center of $G$ acts by a unitary character. If $w_1+w_2$ is even one could also twist by $\det^{-\frac{w_1+w_2}{2}}$.}   
\end{para}

\begin{para}\label{Twisted Translation Functors}{\it Twisted translation functors.} Let $\theta: \Qpx \ra \Ex$ be a locally analytic character. We define the twisted translation functor $\theta \ot \cTclamu$ on the category of coadmissible $\DG$-modules by the formula 

\[(\theta \ot \cTclamu)(M) = (\theta \circ \det) \ot_E\cTlamu(M) \;.\]

\vskip8pt

In the remainder of section \ref{Galois} we assume that 

\begin{numequation}\label{cond theta}
\theta \; \mbox{ satisfies } \;\; v_p(\theta(p)) = -\frac{w_1+w_2}{2} \;.
\end{numequation} 

Such a $\theta$ exists if we pass to a quadratic extension of $E$, if $w_1+w_2$ is odd. In this case  we may take $\theta(x) = \omega^{v_p(x)}$ for a fixed square root $\omega$ of $p^{-w_1-w_2}$.
\end{para}

\vskip8pt

\begin{para}{\it Conditions on $\wnu$.} Set

\begin{numequation}\label{hat defs}
\begin{array}{ccl}\wmu & = &\wlambda \wnu (\theta \ot \theta) = \delta_2^{-1}x^{w_1}\theta \ot \delta_1^{-1}x^{w_2}\theta x|x| \;,\\
&&\\
\wdelta_1 & = & \delta_1 x^{-w_2}\theta^{-1} \;,\\
&&\\
\wdelta_2 & = & \delta_2 x^{-w_1}\theta^{-1} \;,\\
&&\\
\whs & = & (\wdelta_1,\wdelta_2,\infty) \;.
\end{array}
\end{numequation}

Note that $\rmd \wmu = \nu + \lambda$, as in \ref{setup}, and  

\begin{numequation}\label{properties}
\begin{array}{lcl}
\pair{\nu}{\alpha^\vee} & = & \frac{w_1-w_2}{2}\;,\\
&&\\
\pair{\lambda+\rho}{\alpha^\vee} & = & w(s) \;,\\ 
&&\\
\pair{\mu+\rho}{\alpha^\vee} & = & w(\whs)  \;\; = \;\; w(s) + 2\pair{\nu}{\alpha^\vee} \;.
\end{array}
\end{numequation}

In order to ensure that $\whs$ is contained in $\sS_{\irr}$, and is generic, and that the translation functor $\cTlamu$ is an equivalence of categories, we consider the following conditions:

\begin{numequation}\label{cond1}
u(s) + \pair{\nu}{\alpha^\vee} > 0 \;,
\end{numequation}

and

\begin{numequation}\label{cond3}
\mbox{ if } w(s) \in \bbZ_{>0}\,, \mbox{ then }  \; w(s) + \pair{\nu}{\alpha^\vee} - u(s) > 0 \;,
\end{numequation}

as well as 

\begin{numequation}\label{cond2}
\mbox{suppose } w(s) \in \bbZ\,,\mbox{ then: } \mbox{ if } w(s) \left\{\begin{array}{c} >0 \\ =0 \\ <0 \end{array}\right\} \mbox{ then }  w(s) + 2\pair{\nu}{\alpha^\vee} \left\{\begin{array}{c} >0 \\ =0 \\ <0 \end{array}\right\}\;.
\end{numequation}

\vskip8pt

We note that \ref{cond1} and \ref{cond3} imply the first case of \ref{cond2} (when $w(s) \in \bbZ_{>0}$). However, \ref{cond1} and the first case of \ref{cond2} do not imply \ref{cond3} (for example, one may have $u(s) = 2$, $\pair{\nu}{\alpha^\vee} = -1$, and $w(s) = 3$).
\end{para}

\vskip8pt

\begin{prop}\label{translation functors for GL2} Assume \ref{cond2} holds. 

\vskip8pt

(i) Then

\[\theta \ot \cTclamu: \DGcoadlambda \longrightarrow \DGcoadmu\] 
and
\[\theta^{-1} \ot \cTcmula: \DGcoadmu \longrightarrow \DGcoadlambda\] 

\vskip8pt

are equivalences of categories and quasi-inverse to each other. 

\vskip8pt

(ii) There are isomorphisms of $\DG$-modules

\[(\theta \ot \cTclamu)\Big(\DG\ot_{D(B)} E_{\wlambda}\Big) \cong \DG\ot_{D(B)} E_{\wmu}\]

\vskip8pt

and 

\[(\theta \ot \cTclamu)\Big(\DG\ot_{D(B)} E_{w_\circ \cdot \wlambda}\Big) \cong \DG\ot_{D(B)} E_{w_\circ \cdot \wmu} \;.\]
\end{prop}

\vskip8pt

\begin{proof} We see from \ref{properties} that 

\[\lambda \in \Lambda \, \Longleftrightarrow \; w(s) \in \Z \; \Longleftrightarrow  w(\whs) \in \Z \, \Longleftrightarrow \; \mu \in \Lambda \;.\] 

\vskip8pt

and by assumption \ref{cond2} we have

\[\lambda \mbox{ is anti-dominant} \, \Longleftrightarrow \; w(s) \notin \bbZ_{>0} \; \Longleftrightarrow  w(\whs) \notin \bbZ_{>0} \, \Longleftrightarrow \; \mu \mbox{ is anti-dominant} \;.\]

\vskip8pt

If $\lambda$ is not anti-dominant, we can replace $\lambda$ and $\mu$ by $w_\circ \cdot \lambda$ and $w_\circ \cdot \mu$, respectively, both of which are then anti-dominant. Since $w_\circ \cdot \mu - w_\circ\cdot \lambda = w_\circ(\nu)$, the finite-dimensional representations used to define $\cTclamu$ and $\cT^{w_\circ \cdot \mu}_{c,w_\circ\cdot \lambda}$ are the same, and the projection operators are the same, because $\chi_\lambda = \chi_{w_\circ \cdot \lambda}$ and $\chi_\mu = \chi_{w_\circ \cdot \mu}$. Hence $\cTclamu = \cT^{w_\circ \cdot \mu}_{c,w_\circ\cdot \lambda}$. Therefore, we may assume without loss of generality that $\lambda$ and $\mu$ are both anti-dominant. 

\vskip8pt 

Assertion (i) follows from \ref{D-equivalence} once we have seen that conditions (i)-(iii) in \ref{BeGe theorem} are fulfilled. By construction, $\lambda$ and $\mu$ are compatible. And as we just saw, we may assume that both are anti-dominant. The stabilizer of $W^\circ_\lambda$ for the dot-action is non-trivial if and only if the restriction of $\lambda$ to the Cartan subalgebra of diagonal matrices in $\frs\frl_2(F)$ is $-\rho$, which is equivalent to $w(s) = 0$. This implies $w(\whs) = 0$, which in turn means that $\mu$ has stabilizer $W^\circ_\mu = W$. Therefore, the conditions \ref{BeGe theorem} are fulfilled and $\cTclamu$ and $\cTcmula$ are mutually quasi-inverse equivalences, and the same is true for the twisted functors in statement (i).

\vskip8pt

(ii) Suppose first that $\lambda$ and $\mu$ are both anti-dominant. If they are both integral, then $\Wbrla = W_{[\mu]} = W$, and the assertion follows immediately from \ref{translation of induced representations}. If they are both non-integral, then $\Wbrla = W_{[\mu]}$ is trivial, but, as we remarked above, we have $\cTclamu = \cT^{w_\circ \cdot \mu}_{c,w_\circ\cdot \lambda}$ and can therefore apply \ref{translation of induced representations} to $\DG\ot_{D(B)} E_{\wlambda}$ as well as to $\DG\ot_{D(B)} E_{w_\circ \cdot \wlambda}$. Finally, if both $\lambda$ and $\mu$ are not anti-dominant, they are necessarily integral, and we can apply \ref{translation of induced representations} to the functor $\cT^{w_\circ \cdot \mu}_{c,w_\circ\cdot \lambda}$, which is equal to $\cTclamu$, and to $\DG\ot_{D(B)} E_{w_\circ \cdot \wlambda}$ as well as to $\DG\ot_{D(B)} E_{w_\circ \cdot (w_\circ \cdot \wlambda)} = \DG\ot_{D(B)} E_{\wlambda}$.
\end{proof}

\vskip8pt

\begin{theorem}\label{trianguline thm}
Let $s = (\delta_1, \delta_2, \infty)$ be a generic point of $\sS_{\irr}$. Let $\wnu = x^{w_1} \ot x^{w_2}$ be an algebraic character of $\bT$ and $\theta$ a locally analytic character of $\Qpx$ satisfying \ref{cond theta}. Define $\wlambda$ as in \ref{lambdatilde} and $\wmu$, $\whs$ as in \ref{hat defs}. Assume \ref{cond1}, \ref{cond3}, and \ref{cond2} are fulfilled. 

\vskip8pt

\begin{enumerate}
\item[(i)] $\whs$ is a generic point in $\sS_{\irr}$.

\vskip5pt

\item[(ii)] There is an isomorphism of $\DG$-modules

\[(\theta \ot \cTclamu)\Big((\Pi(s)^{\la})'\Big) \cong \Big(\Pi(\whs)^{\la}\Big)' \;.\] 
\end{enumerate}
\end{theorem}

\begin{proof}
(i) We have 

\[\begin{array}{rcl}
v_p(\wdelta_1(p))+ v_p(\wdelta_2(p)) & = & v_p(\delta_1(p) p^{-w_2}\theta(p)^{-1})+ v_p(\delta_2(p)p^{-w_1}\theta(p)^{-1}) \\
&&\\
& = & v_p(\delta_1(p)) - w_1 - v_p(\theta(p))+ v_p(\delta_2(p)) - w_2 - v_p(\theta(p)) \\
&&\\
& = & v_p(\delta_1(p))+ v_p(\delta_2(p)) - (w_1 + w_2) - 2v_p(\theta(p))\\
&&\\     
& = & 0 \;,
\end{array}\]

\vskip8pt

since we assume \ref{cond theta} and because $s \in \sS_*$. Furthermore, we have 

\[u(\whs) = v_p(\delta_1(p)p^{-w_2}\theta(p)^{-1}) = u(s) - w_2 + \frac{w_1+w_2}{2} = u(s) + \pair{\nu}{\alpha^\vee} > 0\]

\vskip8pt

by \ref{cond1}. This shows that $\whs$ is contained in $\sS_*$. Note that 

\[\wdelta_1.\wdelta_2^{-1} = (\delta_1.\delta_2^{-1}) x^{-w_2+w_1} \;.\] 

\vskip8pt

Because we assumed $s$ to be generic, we know that $\delta_1. \delta_2^{-1}$ is not of the form $x^i|x|$ for any $i \in \bbZ_{\ge 1}$, cf. \ref{no *-pts over Z'}. If $\wdelta_1.\wdelta_2^{-1}$ is of the form $x^j|x|$ for some $j \in \Z$, then we must have $\delta_1. \delta_2^{-1} = x^i|x|$ with $i \in \bbZ_{\le 0}$ and $j = i+w_1-w_2$. In this case $w(s) = i$ is an integer $\le 0$. Then condition \ref{cond2} implies that $0 \ge w(\whs) = w(s) + 2\pair{\nu}{\alpha^\vee} = w(s) + w_1-w_2 = i+w_1-w_2 = j$. Therefore $\whs$ is generic. 

\vskip8pt

If $s \in \sS^{\rm ng}_*$, then $w(s) \notin \Z$ or $w(s) \in \bbZ_{<0}$. In the former case $w(\whs) = w(s) + w_1-w_2 \notin \Z$ by \ref{properties}, and in the latter case we have $w(\whs) \le 0$, by \ref{cond2}. Hence $\whs \in \sS_*^{\rm ng}$.

\vskip8pt

Now suppose $w(s) \in \bbZ_{>0}$, hence $s \in \sS_*^{\rm cris}$. Then $w(\whs) \in \bbZ_{>0}$ by \ref{cond2}. Furthermore, by \ref{properties}, $w(\whs) - u(\whs) = w(s) + 2\pair{\nu}{\alpha^\vee} - u(s) - \pair{\nu}{\alpha^\vee} = w(s) + \pair{\nu}{\alpha^\vee} - u(s) >0$, as we assume \ref{cond3}. Therefore, $\whs \in \sS_*^{\rm cris}$. This proves (i).

\vskip8pt

(ii) We apply the translation functor $\theta \ot \cTclamu$ to the exact sequence \ref{colmez 8.7 seq}, and because of \ref{translation exact} we obtain the exact sequence

\begin{numequation}\label{transl on ex seq}
\begin{array}{ccccc}
0 & \ra & (\theta \ot \cTclamu)\Big(\DG\ot_{D(B)} E_{w_\circ \cdot \wlambda}\Big) && \\
&&&&\\
& \ra & (\theta \ot \cTclamu)\Big((\Pi(s)^{\la})'\Big) &&\\
&&&&\\
& \ra & (\theta \ot \cTclamu)\Big(\DG\ot_{D(B)} E_{\wlambda}\Big) & \ra & 0\;. \\
\end{array}
\end{numequation}

As we assume \ref{cond2}, we can apply \ref{translation functors for GL2} and find that \ref{transl on ex seq} gives the exact sequence

\begin{numequation}\label{transl on ex seq 2}
\begin{array}{ccccc}
0 & \ra & \DG\ot_{D(B)} E_{w_\circ \cdot \wmu} && \\
&&&&\\
& \ra & (\theta \ot \cTclamu)\Big((\Pi(s)^{\la})'\Big) &&\\
&&&&\\
& \ra & \DG\ot_{D(B)} E_{\wmu} & \ra & 0\;. \\
\end{array}
\end{numequation}

Passing to strong dual spaces and using \ref{Tensor description of induction} we obtain from \ref{transl on ex seq 2} the exact sequence

\begin{numequation}\label{transl on ex seq 3}
0  \lra B^\la(\wdelta_1,\wdelta_2) \lra  \left[(\theta \ot \cTclamu)\Big((\Pi(s)^{\la})'\Big)\right]'_b \lra B^\la(\wdelta_2,\wdelta_1) \lra  0 \;.
\end{numequation}

By \ref{translation functors for GL2} translation functor $\theta^{-1} \ot \cTcmula$ is quasi-inverse to $\theta \ot \cTclamu$, the sequence \ref{transl on ex seq 3} is also not split. By \ref{colmez 8.7} (ii), this means that the locally analytic representation in the middle is unique up to isomorphism, cf. also \cite[8.8. (i)]{Co2}. Therefore,

\[\left[(\theta \ot \cTclamu)\Big((\Pi(s)^{\la})'\Big)\right]'_b \; \cong \; \Pi(\whs)^\la \;\]

\vskip8pt

which after passing to strong dual spaces gives the assertion.
\end{proof}

\begin{rem} It seems natural to expect that a result similar to \ref{trianguline thm} also holds for trianguline representations $V(s)$ where $s$ special (and not generic, as assumed above). Furthermore, when $V$ is de Rham with irreducible associated Weil group representation, the effect of the translation functors on $\Pi(V)^\la$ (or rather its dual space) can be interpreted in the framework of the change of weights theory of Colmez in \cite{Colmez_poids}.  \end{rem}

\bibliographystyle{abbrv}
\bibliography{mybib}

\end{document}